\documentclass[11pt]{article}
\usepackage[margin=1.2in,footskip=0.25in]{geometry} 
\usepackage{lineno,hyperref}
\modulolinenumbers[5]

\usepackage{endnotes}
\let\footnote=\endnote

\usepackage{algorithm,algpseudocode}
\usepackage{mathrsfs}
\usepackage{xcolor}
\usepackage{longtable}
\usepackage{url} 
\usepackage{footnote}
\usepackage{lscape} 
\usepackage{threeparttablex}
\usepackage{graphicx,subfigure}
\usepackage{amsmath}
\usepackage{amsthm}  
\usepackage{mathtools}
\usepackage{amssymb}
\usepackage{caption}
\usepackage{accents}
\usepackage{algorithm,algpseudocode}
\usepackage{booktabs}
\usepackage{comment}

\newtheorem{definition}{Definition}
\newtheorem{theorem}{Theorem}
\newtheorem{proposition}{Proposition}
\newtheorem{lemma}{Lemma}
\newtheorem{corollary}{Corollary}
\newtheorem{remark}{Remark}

\numberwithin{observation}{section}
\numberwithin{assumption}{section}
\numberwithin{definition}{section}
\numberwithin{theorem}{section}
\numberwithin{proposition}{section}
\numberwithin{lemma}{section}
\numberwithin{corollary}{section}

\DeclarePairedDelimiterX\Set[2]{\lbrace}{\rbrace}%
 { #1 \,\delimsize| \,\mathopen{} #2 }

\newcommand{\E}{\mathbb{E}}
\newcommand{\mtc}{\mathcal}

\newcommand{\W}{\mathcal{W}}
\newcommand{\bs}[1]{\boldsymbol{#1}}

\newcommand\norm[1]{\lVert#1\rVert}

\newcommand{\beq}{\begin{equation}}
\newcommand{\eeq}{\end{equation}}
\newcommand{\ba}{\begin{aligned}}
\newcommand{\ea}{\end{aligned}}
\newcommand{\bdm}{\begin{displaymath}}
\newcommand{\edm}{\end{displaymath}}

\newcommand{\mP}{\mathcal{P}}

\newcommand{\conv}{\textrm{conv}}
\newcommand{\mF}{\mathcal{F}}

\newcommand{\mS}{\mathcal{S}}

\newcommand{\whb}[1]{\widehat{\bs{#1}}}

\newcommand{\hbeta}{\hat{\bs{\beta}}}
\newcommand{\Puy}{\mathcal{P}^{\bs{y}}_{u}}
\newcommand{\Puyij}{\mathcal{P}^{\bs{y}}_{u_{ij}}}

\newcommand{\bvij}{\bs{v}^{ij}}
\newcommand{\bvijo}{\bs{v}^{ij1}}
\newcommand{\bvijt}{\bs{v}^{ij2}}

\newcommand{\Uij}{U^{ij}}
\newcommand{\Uijo}{U^{ij1}}
\newcommand{\Uijt}{U^{ij2}}
\newcommand{\Rij}{R^{ij}}
\newcommand{\sij}{s^{ij}}

\newcommand{\cvd}{COVID-19 }

	\newcommand{\blind}{0}
	
    \makeatletter
    \renewcommand\section{\@startsection {section}{1}{\z@}%
                                       {-3.5ex \@plus -1ex \@minus -.2ex}%
                                       {2.3ex \@plus.2ex}%
                                       {\normalfont\fontfamily{phv}\fontsize{16}{19}\bfseries}}
    \renewcommand\subsection{\@startsection{subsection}{2}{\z@}%
                                         {-3.25ex\@plus -1ex \@minus -.2ex}%
                                         {1.5ex \@plus .2ex}%
                                         {\normalfont\fontfamily{phv}\fontsize{14}{17}\bfseries}}
    \renewcommand\subsubsection{\@startsection{subsubsection}{3}{\z@}%
                                        {-3.25ex\@plus -1ex \@minus -.2ex}%
                                         {1.5ex \@plus .2ex}%
                                         {\normalfont\normalsize\fontfamily{phv}\fontsize{14}{17}\selectfont}}
    \makeatother
	
	\usepackage{enumerate}
	\usepackage{natbib} 
	
\begin{document}
		
		\def\spacingset#1{\renewcommand{\baselinestretch}%
			{#1}\small\normalsize} \spacingset{1}
		
		\if0\blind
		{
			\title{Service Center Location Problem with Decision Dependent Utilities with an Application to Early Stage Testing and Vaccination in Epidemic Planning}
			\author{Fengqiao Luo and Sanjay Mehrotra \\
			Department of Industrial Engineering and Management Science, \\
			Northwestern University, Evanston, IL, 60202 }
			\date{}
			\maketitle
		} \fi
		
		\if1\blind
		{

            \title{Service Center Location Problem with Decision Dependent Utilities}
			\author{Author information is purposely removed for double-blind review}
			
\bigskip
			\bigskip
			\bigskip
			\begin{center}
				{\LARGE\bf \emph{IISE Transactions} \LaTeX \ Template}
			\end{center}
			\medskip
		} \fi
		\bigskip
		
	\begin{abstract}
We study a service center location problem with ambiguous utility gains upon receiving service. The model is motivated by the problem of deciding medical clinic/service centers, possibly in rural communities, where residents need to visit the clinics to receive health services. A resident gains his utility based on travel distance, waiting time, and service features of the facility that depend on the clinic location. The elicited location-dependent utilities are assumed to be ambiguously described by an expected value and variance constraint. We show that despite a non-convex nonlinearity, given by a constraint specified by a maximum of two second-order conic functions, the model admits a mixed 0-1 second-order cone (MISOCP) formulation. We study the non-convex substructure of the problem, and present methods for developing its strengthened formulations by using valid tangent inequalities. Computational study shows the effectiveness of solving the strengthened formulations. Examples are used to illustrate the importance of including decision dependent ambiguity. An illustrative example to identify locations for Covid-19 testing and vaccination is used to further illustrate the model and its properties. 
	\end{abstract}
			
	\noindent%
	{\it Keywords:} Facility location; Robust optimization; Utility Functions.

	\spacingset{1.5} 

\section{Introduction}
\label{sec:introd} 
This paper considers a location problem that decides to locate service centers among candidate locations to serve customers from different sites. The objective is to maximize the total utility of the service received by the customers. The utility gained by the customers depends on the \textit{location pattern} of the service centers. The decision is constrained due to available budget, and other considerations such as staff availability. More formally, let $S$ be the index set of customer sites and $F$ be the index set of candidate service center locations. Let $h_j$, $j\in F$ be the cost of opening a service center, and $y_j$ be the binary decision variable for opening a service center at location $j\in F$. The customer site $i\in S$ has a demand $D_i$. Each service center $j\in F$ has a service capacity $C_j$. $x_{ij}$  represents the coverage of demand generated from customer site $i$ to facility location $j$.  We let $u_{ij}(\bs{y})$ represent the utility gained by a customer at site $i\in S$ receiving service from the service center at location $j\in F$, if service center locations are given by $\bs{y}$, where $\bs{y}=\Set*{y_j}{j\in F}$.  The utility gain is ambiguous, and we let $\Puyij$ represent the ambiguity set of utility functions $u_{ij}(\bs{y})$ for $i\in S$, $j\in F$ when the location decision $\bs{y}$ is implemented. Note that this utility gain, as well as the ambiguity set describing the utility gain, is dependent on the decision vector $\bs{y}$.
The model assumes a total budget $B$ for the service center decisions. For a general model we assume that $h_j$ is an additional fixed gain for opening a service center $j\in F$. We can take $h_j = 0$ for a problem that only needs to decide the location of service centers. The utility-robust service center location model formulation is given as follows:
\begin{equation}\label{opt:RFL}
\begin{aligned}
&\underset{\bs{y}}{\textrm{max}}\;\;\bs{h}^{\top}\bs{y}+\mathcal{Q}(\bs{y}) \\
&\textrm{ s.t. } \sum_{j\in F} b_j y_j \le B,   \\
&\qquad y_j\in\{0,1\}\;\forall j\in F,
\end{aligned}
\tag{RFL}
\end{equation}
where the constraint in \eqref{opt:RFL} is the budget constraint. Additional structural constraints on $\bs{y}$ may be included though they are not given here.  For a location decision vector $\bs{y}$, 
$\mathcal{Q}(\bs{y})$ is a risk-averse utility gain given by the following problem (RSP):
\begin{equation}\label{opt:RFL-II}
\begin{aligned}
\mathcal{Q}(\bs{y})=\; & \underset{\bs{x}}{\textrm{max}}\;\sum_{i\in S}\sum_{j\in F}
x_{ij}\left(\underset{P\in\Puyij}{\textrm{min}}\;\E_P[u_{ij}(\bs{y})] \right) \\
&\textrm{ s.t. } \bs{x}\in X(\bs{y}),
\end{aligned}
\tag{RSP}
\end{equation}
where the feasible set $X(\bs{y})$ is defined as 
\begin{equation}\label{eqn:X-set}
X(\bs{y}):=\Set*{x_{ij}\;\forall i\in S,\;\forall j\in F}
{\def\arraystretch{1.2}
\begin{array}{ll}
\sum_{i\in S}x_{ij} \le C_jy_j & \forall j\in F, \\
\sum_{j\in F}x_{ij} \le D_i & \forall i\in S, \\
x_{ij}\ge 0 & \forall i\in S,\; \forall j\in F
\end{array}
}.
\end{equation}
The objective of \eqref{opt:RFL-II} is to maximize the worst-case expected
\textit{maximum potential utility} gained. To estimate the maximum 
potential utility under a given location decision, it is assumed customers from all
sites will collaboratively share the limited capacity from opened facilities to maximize
the total utility gain. In this sense, the customer flows $\bs{x}$ are
treated as virtual variables that can be determined by a central policy maker
in the model. Despite very ideal, this assumption is valid if the goal
is to estimate the limit of utility gain led by a given location pattern. 
Without the assumption, one needs to specify a lot more mechanism under limited capacity of service,
such as priority of providing service, a fare policy of sharing capacity among
customers from different sites, and the relation between the customer flow 
and the utility etc., which is beyond the scope of this paper.

The source of ambiguity in \eqref{opt:RFL} is from the evaluation of expected utility $u_{ij}(\bs{y})$.
The first constraint in \eqref{eqn:X-set} is the capacity constraint for each service center, and the second constraint ensures that the total number of customers from site $i$ cannot exceed the potential demand from $i$. Note that the results in this paper remain valid when the set $X(\bs{y})$ is defined differently from an alternative application.

\subsection{Possible applications of the modeling framework}
\label{sec:model-appl}
The model studied in this paper is motivated by the situations where customers 
go to a service center in order to receive service.
The utility of service received by customers is effected by the joint locations of the service centers. This feature makes \eqref{opt:RFL}  different from the traditional facility location problems in which resources are delivered from a facility to customers
to meet demand, and a delivery cost is incurred \citep{book_daskin_disc-loc}.
We give some real world situations to which our model can be applied. 

In the first example, we consider a healthcare system of a developing country where the state and central governments plan to open primary care clinics with a limited budget e.g., as in \citep{sharma2016-commun-clinic}. 
The clinics provide primary care and health screening for the residents at low or no cost. Since patients need to come to a clinic to receive healthcare services, the value of these clinics to a resident depends on the location of the clinics, 
especially the distance and accessibility from the place of residence. Residents have a choice of clinic, and may go to multiple clinics. Each clinic has a limited capacity. As discussed above, the utility of the clinics to a resident depends on their locations.
Analogous examples arise in the context of opening low cost subsidized pharmacies, fitness centers, or testing locations. The latter is used as a case study in the context of deciding Covid-19 vaccination locations. 

In the context of for-profit organizations, consider the problem of locating a few shopping centers in a city.
Different locations and features (i.e., scale, presentation, neighborhood and quality of service) 
of shopping centers may attract the residents differently, which results in a location dependent shopping center experience (utility) gain. Since merchandise selling price is typically matched, it  is in not necessarily the primary difference of the shopping centers from it competitors.

\subsection{Contributions of this paper}
This paper makes the following contributions:
\begin{itemize}
	\item We establish a utility-robust optimization model \eqref{opt:RFL} for the service center location problem when 
		  the utilities are random parameters with ambiguous probability distribution that
		  are affected by the service center locations. Under a suitable moment-based model for specifying location
		  utility ambiguity set, we show that it is possible to reformulate \eqref{opt:RFL} as 
		  a mixed 0-1 second-order cone program (MISOCP).
	\item We investigate the properties of the non-convex constraint, written as the max of two second-order-cone functions,
		  arising in the reformulation of the ambiguity set for the utilities. We give representations
		  of the convex hull associated with the non-convex constraint.  
	\item We develop numerical frameworks for generating tangent inequalities of the convex hull associated with the non-convex constraint. 
		  These tangent inequalities lead to stronger formulations of \eqref{opt:RFL}.
	A numerical study is conducted to test the computational performance of solving \eqref{opt:RFL} instances 
	       with or without the convexification cuts developed in this paper. 
	       Computational results show that incorporating convexification cuts results in a significant cpu time savings. It allows us to solve problems with up to 3,000 potential locations and 300 site budget in less than 1/2 hour. 
	       
	       \item Numerical experiments are used to illustrate properties of the \eqref{opt:RFL} model and discuss insights. An illustrative example to identify Covid-19 test center locations is used to further illustrate the model and its properties. 
\end{itemize}

\subsection{Organization of this paper}
This paper is organized as follows. Section~\ref{sec:liter-rev} provides a literature review on the
facility location problems.
Section~\ref{sec:mod-interp} provides a rationale for the utility's dependence on the
service center locations. Section~\ref{sec:mod-ass} discusses a linear utility assumption we use to
model the decision dependent utility in this paper. Section~\ref{sec:momt-amb-set} establishes an ambiguity set of utilities based on the first two moments of the random utility function.
Section~\ref{sec:illus-exp} provides an illustrative example to show that the robust optimal
service center locations can change with respect to different ambiguity level in the utility.

Section~\ref{sec:mixed-bin-conic} presents a mixed 0-1 second-order-cone program (MISOCP) reformulation of \eqref{opt:RFL}.
In Section~\ref{sec:gen-conv-hull}, we investigate the properties of 
the non-convex constraint in the formulation and give two representations
of the convex hull associated with the non-convex constraint.
Based on the two representations of the convex hull, 
we develop numerical methods for generating tangent inequalities of the convex hull associated with the non-convex constraint.
In Section~\ref{sec:comp-determ-dmd}, 
we provide our computational experience with the MISOCP reformulation of \eqref{opt:RFL} 
and the effectiveness of cuts developed in this paper for solving 41 \eqref{opt:RFL} instances
ranging from small size to large size.  
In Section~\ref{sec:comp-manag-ins}, we discuss some insights in the optimal location changes as a consequence of the utility ambiguity levels. This is followed by the concluding remarks section, where we present a generalization of the model that allows random demand.

\subsection{Literature review}
\label{sec:liter-rev}
\subsubsection{Facility location models and endogenous ambiguity}
Facility location models are extensively studied in \citep{book_daskin_disc-loc}. In the facility location problem, a decision maker needs to decide location of a limited number of  facilities (factories, retail centers, power plants, service centers, etc.), and determine coverage of demand from different sites by the located facilities. The objective is to minimize the facility setup cost and the cost of production/delivery. The facility location models provide framework for other problems in resource allocation, 
supply chain management and logistics, etc. \citep{2009-fac-loc-rev}.

\cite{2003-rob-fac-loc} investigated the problem of locating a single facility in a continuous region of $\mathbb{R}^2$ that meets the demand. The facility is located in a robust sense by selecting a location that minimizes the perturbed delivery cost with respect to a reference demand distribution. 
\cite{2010-fac-loc-rob-opt-approach} studied a robust multi-period facility location problem with a   box uncertainty set and an ellipsoidal uncertainty set of demand in each period. The model is reformulated as a mixed 0-1 linear program and a mixed 0-1 conic quadratic program, respectively. The objective is to maximize the total profit.  The numerical study showed that robust models provide small but significant improvements over the solution to the deterministic model using nominal demand. 
\cite{2014-rob-fac-loc-hazard-waste-transport} analyzed a robust hazardous material carrier allocation problem with a box uncertainty set for the amount of hazardous material in a finite set of sites and for the exposure risk at each link during transport. Here the objective is to minimize the weighted combination of the facility opening/setup cost, delivery cost and total risk exposure. The problem is reformulated as a mixed 0-1 linear program using linearization techniques.  

In stochastic programming based facility location models, the uncertain demand is modeled as a discrete random variable on a finite set of scenarios. Specifically, \cite{1992-dual-stoch-fac-loc} provided an early investigation on a two-stage stochastic optimization model of the uncapacitated
facility location problem with recourse when demand, selling price, production and transportation costs are random. 
\cite{2002-alg-fac-loc-stoch-demand} developed an immobile server location model which is motivated from the problem of locating bank ATMs
or Internet mirror sites congested by stochastic demand originating from nearby customer locations. Here the queueing system for 
each server is modeled by an M/M/1 queue, and the objective is to minimize customers' total travel and waiting time.
\cite{2006-reg-stoch-fac-loc} proposed an $\alpha$-reliable mean-excess regret model ($\alpha$-RMERM) for stochastic facility location modeling.
For a decision $\bs{y}$ of facility location, the regret under a scenario is defined as the increased value in the total weighted delivery distance under 
the decision $\bs{y}$ compared to the minimum value under scenario $s$.
In comparison with the previous $\alpha$-reliable minimax model ($\alpha$-RMM) that minimizes the $\alpha$ quantile of regrets,
the $\alpha$-RMERM minimizes the expectation of the excess regret with respect to the $\alpha$ quantile. Since the mixed 0-1 programming reformulation of the $\alpha$-RMERM is more compact (no big-M coefficient) than that of the $\alpha$-RMM model, 
$\alpha$-RMERM is shown to be computationally more efficient.  
A two-stage stochastic facility location model is also developed for humanitarian relief logistics \citep{2012-two-echelon-stoch-fac-loc}
to minimize the total cost of rescue center location, inventory holding, transportation and shortage of relief items. 

In recent years, research on robust and stochastic facility location (RSFL) models has
investigated a supply chain network where each demand site is allowed to source
supply from multiple distribution centers
\citep{2017-multi-source-suppl-chance-constr}. The objective is to minimize the total cost while satisfying the demand with a given probability. \cite{2017-multi-source-suppl-chance-constr} proposed a set-wise approximation and reformulated the chance constraint in the model using exponentially many second-order cone constraints. A mixed binary second-order cone program is solved numerically using a cutting plane procedure. 
\cite{2017-rob-deploy-card-arrest-loc-uncert} studied a distributionally robust medical equipment (defibrillators) location problem to reduce cardiopulmonary resuscitation (CPR)
delay in sudden cardiac arrest patients due to the defibrillator distance from the event site. Based on the defibrillator location, 
the objective of this model uses conditional value-at-risk (CVaR) on the distance between 
the cardiac arrest event site and the nearest defibrillator location. The uncertainty set of the cardiac arrest event site is constructed using a finite set of possible locations. 
It is shown that this model can be reformulated as a mixed 0-1 semi-infinite program, and a row-and-column generation algorithm is applied to solve the reformulated problem.   

A recent trend in robust and distributionally-robust optimization is to incorporate the endogenous
ambiguity in the modeling framework, as in many real-world decision-making systems the uncertainty
of a system is likely to be decision-dependent. \cite{nohadani2016_opt-dec-uncert} investigated
the the reformulation of robust linear programs (possibly involving discrete variables) 
with polyhedral decision-dependent uncertainty set. 
This modeling framework has been applied to investigate a robust shortest-path 
problem in which selection of road links have impact 
on the size and shape of the following-up uncertainty set.  
\cite{luo-2018-D3RO} and \cite{noyan2018} have investigated
distributionally-robust optimization models with decision-dependent
ambiguity set for the candidate unknown probability distributions of
model parameters. On the application side, the endogenous uncertainty
has been handled in resource management \citep{tsur2004_gdwater-threat-catastr-event}, 
stochastic traffic assignment modeling \citep{shao2006_reliab-stoch-traff-assign-demand-uncert}, 
oil (natural gas) exploration \citep{jonsbraten-thesis1998_opt-mod-petrol-field-expl,
grossmann2009_stoch-prog-plan-oil-infrast-dec-uncert,
grossmann2004_stoch-prog-plan-gas-field-uncert-reserv}, 
and robust network design \citep{ahmed2000-phd_plan-uncert-stoch-MIP,
viswanath2004_invest-stoch-netwk-min-exp-short-path}.

The endogenous ambiguity for distributionally-robust optimization
has been investigate with applications on facility location problems.
\cite{basciftci2021} investigated a decision-dependent distributionally-robust 
(D3RO) facility location problem, in which the ambiguity set is defined 
using disjoint lower and upper bounds on
the mean and variance of the candidate probability distributions for customer demand,
and these bounds are defined as linear functions of the location vector to admit
a tractable MILP reformulation. \cite{luo2020} investigated the decision-dependent
customer demand from a preference-level point of view. This approach admits a 
decoupled way of incorporating decision dependence and distributional robustness into the model,
which significantly improves the computational efficiency. The problem investigated
in this paper has a different setting from \citep{basciftci2021} and \citep{luo2020}.
The model in \citep{basciftci2021} is established on a classic facility location problem in which
resources are delivered to multiple customer locations to meet the demand with decision-dependent
ambiguity. In contrast to \citep{basciftci2021}, the model in this paper is established for the case that
customers are required to visit the facilities in order to get service. The quality
of service and demand fulfillment are both taken into account in the model. 
The model constructed in \citep{luo2020} is based on a so called ``maximum attraction principle''
assumption to establish a more careful way of counting the number of customers 
who are willing to visit an opened facility, but this assumption needs a further justification
with evidences from a real-world example. The development of this paper is not based on 
that assumption.

\subsubsection{Decision theory and utility models}
Utility models are widely used in economics and consumer theory for decision making based on discrete choices
\citep{fishburn1970_utility-thy-dec-making,dyer1992_multi-criteria-dec-making-utility,zavadskas2011_multi-criteria-dec-making-rev}.
\cite{fishburn1970_utility-thy-dec-making} provided a fundamental understanding of utility theory for decision making,
focusing on the logic of utility comparison and the structure of utility functions.
There are several classes of utility modeling frameworks, among which the expected utility models 
 is commonly used \citep{schoemaker1982_exp-utility-model}. 
The expected utility theory is based on assumptions including independent evaluations,
exhaustive search, trade-offs, objective probabilities and values,
which helps simplify the modeling of a complex psychological process of decision making \citep{katsikopoulos2008_one-reason-dec-making}.  
An expected utility model evaluates multiple choices based on some attributes. Every choice has a value
at each attribute, and the decision maker is characterized by a weight vector (independent of choices) of the attributes  
which describes the preference levels of the attributes. The choice that maximizes the expected utility is used as the optimal choice.

\cite{luce1991_linear-utility-models-bin-gambles} studied linear utility models for binary decision making.
\cite{bell1982_reg-dec-making-uncert} incorporated regret into utility function and used numerical examples to show that
this modification can improve prediction and lead to better description of decision makers' behavior in some situations. 
\cite{rabin1997_risk-aversion-exp-utility-thy} studied a preliminary risk-averse random utility model. \cite{friedman2003_learn-prob-model-exp-utility-max} investigated a problem of learning a probabilistic model 
based on maximizing the expected utility with prior knowledge on the unknown probability of events.
\cite{cascetta2009_dom-among-alter-rand-utility-models} extended the concept of dominance among
alternatives to the framework of the random utility theory.
\cite{katsikopoulos2008_one-reason-dec-making} analyzed a decision-making utility model  
for some realistic situations where only a few attributes play a dominant role, 
and exhaustive computation is not achievable for the decision maker.
\cite{kitamura2018_nonparam-analy-rand-utility-md} developed and implemented a nonparametric test of random utility models.
\cite{huang2013_agg-utility-group-dec-making} investigated an approach for group decision making that is based on
aggregating individual utility models. 
The linear utility model is a building block for establishing several probabilistic choice models \citep{mcfadden2000_mixed-mnl-model}.

Utility models are widely applied in research areas such as social choice \citep{soufani2012_rand-utility-thy-social-choice},
welfare analysis and comparison among individuals \citep{decoster2010_rand-utility-model-labor-supply},
nursing practice measurement \citep{brennan2000_nursing-pract-utility}, 
simulation and estimation of travel demand, travel time and route choice optimization
\citep{cascetta2001_rand-utility-simul-travel-dmd,blayac2001_utility-travel-time,hawas2004_route-choice-utility-models}, 
early drug discovery \citep{parrott2005_utility-early-drug-discovery}, and sustainable forest management 
\citep{wintle2005_utility-dyn-landscape-sust-forest-managmt}, etc.
Nondecreasing concave utility function is used as an objective value in risk-averse decision making in financial market models 
\citep{rasonyi2005_utility-max-discrete-time-finan-market-models,mehrotra2015-rob-portf-opt-utility}.

The use of utility in more complex decision-making models has also received some attention in the robust optimization literature. 
\cite{ahmed2000-phd_plan-uncert-stoch-MIP} investigated a class of single-stage stochastic programs 
with discrete candidate probability distributions that are based on Luce's choice axiom~\citep{luce1977-choice-axiom}. 
\cite{schied2005_opt-invest-rob-utility} studied an optimal investment strategy based on a distributionally-robust utility model.
\cite{mehrotra2015-rob-portf-opt-utility} studied a model that searches for a robust optimal decision over a set of risk-averse utilities.  
This modeling framework is further extended to the context of general utilities in \citep{hu2018-rob-dec-making-gen-utility}.

\subsubsection{Difference with some of the related works}
Different versions of facility location problems with customer utilities have been investigated in the existing literature. For example, the early work of \cite{benati2002} investigated a problem of allocating new facilities from a newcomer which enters a market and will compete with the existing competitors. Utility of customers is modeled using a multinomial logit (MNL) based choice model, which enters an integer programming to find an optimal decision of locations. A similar setting is investigated in \citep{ljubic2018}, with some submodular cuts being derived to strengthen the formulation and improve the computational performance. The MNL choice model has also been used in preventative healthcare facility location models \citep{haase2015} and healthcare network design problems \citep{denoyel2017}. The work in \citep{garcia2018} investigated a robust location problem of new housing developments using the MNL choice model. It focused the impact of the uncertainty in demand and customer utilities on the allocation of new houses. It handled the two source of uncertainty by introducing protection terms defined based on a budget uncertainty set to incorporate them in the constraint and objective.

There are a few key differences in this paper compared with the above literature.
First, the total utilities in our model are the objective to be optimized, whereas they  
are used as arguments in the MNL choice model to determine the customer flow 
in the above literature. Furthermore, our goal in this paper is not to find facility
locations with a realistic model of customer flow determined by the choice model, 
but to find facility locations that lead to the maximum potential total utility gain from
customers. This logic is discussed in Section~\ref{sec:mod-interp} with more details.

\section{Decision Dependence in Facility Location Utility Assessment}
\label{sec:DRO-FL}
\subsection{Model Interpretation}
\label{sec:mod-interp}

To give an example that the utility gain from customers getting service can potentially depend
on the location pattern of certain facilities, we consider a situation in which there is only one opened service center in total and it is located at the place $j\in F$. The utility for customers at site $i$ of getting service from facility $j$ in this case is denoted by $u_{ij}(\bs{e}_j)$ according to the model, where $\bs{e}_j$ is a $|F|$-dimensional vector with the $j^{\textrm{th}}$ entry being one and other entries being zero. The utility $u_{ij}(\bs{e}_j)$ can be overestimated or underestimated, 
since no other service centers are available for comparison. If a second service center is opened at location $j^{\prime}$, customers can compare the traveling distance, level and quality of service from the two service centers ($j$ and $j^{\prime}$). Based on this comparison, customers at site $i$ may modify their utility value for $j$. In other words, the \textit{subjective} utility viewed by customers from site $i$ towards getting service at the facility $j$ could be impacted by the presence of $j^\prime$ especially in the situation that $j$ and $j^\prime$ are located in the neighborhood of $i$. Therefore, in general the utility $u_{ij}$ can depend on the 
location pattern $\bs{y}$ of opened facilities. In practice, the utility gain is a very subjective quantity that relies on non-rigorous comparison among opened facilities by customers, which justifies the subtle dependency on the locations of facilities. Another concrete example is about the restaurant choice in a small town. If there is only one restaurant in the town. Customers are more likely to give a high rate to it, as they have no other options to compare. But if there are multiple ones, customers could be more selective to score those restaurants. 
   
To make the evaluation of utility more realistic, the location-pattern dependent 
utility $u_{ij}(\bs{y})$ can be modeled as a random function, 
and the probability measure of this random function depends on the customer site $i\in S$,
the service center $j\in F$, and the location decision vector $\bs{y}$.
Given a probability measure $P\in\Puyij$ of the utility $u_{ij}(\bs{y})$, 
the expected value of $u_{ij}(\bs{y})$ is evaluated from
\begin{equation}
\E_{u_{ij}\sim P}[u_{ij}] = \int_{\mathbb{R}_+} u P(du),
\end{equation}
where $P$ is an element of the ambiguity set $\Puyij$. Note that in the
distributionally-robust setting, it is the ambiguity set (i.e., the size and shape) that depends on $\bs{y}$.
But an element in the set, which is a probability measure, does not explicitly depend on $\bs{y}$.

\subsection{Model Assumptions}
\label{sec:mod-ass}
 We assume that $u_{ij}(\bs{y})\ge 0$, and  $u_{ij}(\bs{y})$ is linear in $\bs{y}$, i.e., 
 $$u_{ij}(\bs{y})=(\bs{\beta}^{ij})^{\top}\bs{y}+\varepsilon^{ij},$$
where $\bs{\beta}^{ij}$ is a $|F|$ dimensional random vector and $\varepsilon^{ij}$ is an error term. 
This assumption simplifies our presentation, though the modeling framework allows for the use of a more general functional form.  

In practice the coefficient $\bs{\beta}^{ij}$ in the utility function description may be estimated using a randomized design, or some other alternative methodology. A randomized design is described below. For each $i\in S$ and $j\in F$, 
we randomly select $N$ residents from site $i$ and generate $N$ random location decision vectors 
$\{\bs{y}^k\}^N_{k=1}$ satisfying $y^k_j=1\;\;\forall k\in [N]$. For the $k^{\textrm{th}}$
selected resident, we ask the resident to score, in the range from 0 to 100, the utility of being assigned to service center $j$ for the location vector $\bs{y}^k$.   
Suppose the score given by the $k^{\textrm{th}}$ selected resident is $s^k$, $k\in [N]$.
Then we can estimate the coefficient vector $\bs{\beta}^{ij}$ using the following linear regression model
\begin{equation}\label{eqn:linear-reg}
[\bs{y}^1,\bs{y}^2,\ldots,\bs{y}^N]^{\top}\bs{\beta}^{ij} = [s^1,s^2,\ldots,s^N]^{\top}+ \bs{\epsilon}^{ij},
\end{equation}
where $\bs{\epsilon}^{ij}$ is an error vector.

\subsection{A Utility Ambiguity Set}
\label{sec:momt-amb-set}
The linear model $\hat{u}_{ij}(\bs{y})=(\bs{\beta}^{ij})^{\top}\bs{y}$ is an estimation of the unknown true utility.
The linear utility model is ambiguous due to uncertainty in  $\bs{\beta}^{ij}$, for $i\in S$ and $j\in F$, possibly due to response 
bias, insufficient sampling, and the choice of linear model (model misspecification). Thus we may be interested in robustifying against the ambiguity in the estimation of $u_{ij}(\bs{y})$. 
Below we present an approach to construct an ambiguity set $\Puyij$ based on the mean vector and covariance matrix 
of estimated $\bs{\beta}^{ij}$. Since this approach is identical for every $i\in S$ and $j\in F$, we omit the indices $i,j$ to simplify the notation.

In our approach we treat $\bs{\beta}$ as a random vector that follows an unknown
probability measure. The linear regression model provides a reference mean vector $\hat{\bs{\beta}}$ and a reference
covariance matrix $\widehat{\bs{\Sigma}}$ of $\bs{\beta}$. 
Let $\hat{\bs{\beta}}^*$ and $\widehat{\bs{\Sigma}}^*$ be the true mean
vector and the true covariance matrix of $\bs{\beta}$. 
Suppose we have an uncertainty set $\mathcal{B}$ of $\hat{\bs{\beta}}^*$ and 
an uncertainty set $\mathcal{E}$ of $\widehat{\bs{\Sigma}}^*$, satisfying
$\hat{\bs{\beta}},\hat{\bs{\beta}}^*\in\mathcal{B}$
and $\widehat{\bs{\Sigma}},\widehat{\bs{\Sigma}}^*\in\mathcal{E}$. 
We now define the ambiguity set $\Puy$ as follows:
\begin{equation}\label{eqn:amb-P}
\Puy=\Set*{P\in (\Xi,\mathcal{F})}
{\begin{aligned}
&\underset{\bs{\beta}\in\mathcal{B}}{\textrm{min}}\;\bs{\beta}^{\top}\bs{y} 
\le \E_{u\sim P}[u] \le \underset{\bs{\beta}\in\mathcal{B}}{\textrm{max}}\;\bs{\beta}^{\top}\bs{y} \\
&\underset{\bs{\Sigma}\in\mathcal{E}}{\textrm{min}}\;\bs{y}^{\top}\bs{\Sigma}\bs{y} \le \E_{u\sim P}\big[\big(u-\bs{\hat{\bs{\beta}}^{\top}\bs{y}}\big)^2\big] \le 
\underset{\bs{\Sigma}\in\mathcal{E}}{\textrm{max}}\;\bs{y}^{\top}\bs{\Sigma}\bs{y}
\end{aligned} }.
\end{equation}
The above ambiguity set restricts the candidate  mean and empirical variance of the utility within a confidence region. Specifically, for a candidate probability measure $P$ the quantity $\E_{u\sim P}[u]$ is the mean of the utility if $u$ follows
the probability measure $P$. The quantity $\E_{u\sim P}\big[\big(u-\bs{\hat{\bs{\beta}}^{\top}\bs{y}}\big)^2\big]$ is interpreted as the variance of $u$ with the mean value estimated using the reference mean.  
The first inequality in \eqref{eqn:amb-P} imposes that for any candidate probability measure $P$,
the mean (based on $P$) of the utility should be upper and lower bounded by the maximum and the minimum  
values of $\bs{\beta}^{\top}\bs{y}$ respectively over the choice of coefficients $\bs{\beta}$ from the uncertainty set $\mathcal{B}$. 
Similarly, the second inequality imposes an upper and a lower bound on the empirical variance of the utility over the  covariance matrix $\bs{\Sigma}$ from the uncertainty set $\mathcal{E}$. We note that the specification of the set $\Puy$ may be generalized to consider bounds based on higher order moment considerations (see e.g., \citep{mehrotra2015}). However, solution of models based on such a definition is beyond the scope of the current paper.

The definition of the ambiguity set in \eqref{eqn:amb-P} is independent of the choice of uncertainty sets $\mathcal{B}$ and $\mathcal{E}$. 
We now propose a specific $\mathcal{B}$ and $\mathcal{E}$ for use in \eqref{eqn:amb-P}. 
We define $\mathcal{B}$ as an ellipsoid set with the center at $\hbeta$,
and define $\mtc{E}$ as set of positive semi-definite matrices with lower and upper bounds obtained based on $\bs{\widehat{\Sigma}}$:
\begin{equation}\label{def:B&E}
\begin{aligned}
&\mathcal{B}=\{\bs{\beta}\in\mathbb{R}^{|F|}\; |\; (\bs{\beta}-\hbeta)^{\top}\bs{A}(\bs{\beta}-\hbeta)\le b^2 \}, \\
&\mathcal{E}=\{\bs{\Sigma}\in\mathbb{R}^{|F|\times |F|}\;|\;\gamma_1\widehat{\bs{\Sigma}} \preccurlyeq \bs{\Sigma}\preccurlyeq \gamma_2 \widehat{\bs{\Sigma}} \},
\end{aligned}
\end{equation}
where $\bs{A}$ is positive semi-definite matrix, and $b$ is a positive parameter. The matrix $\bs{\Sigma}$ is positive semi-definite, and scalar $\gamma_1,\gamma_2$ are positive parameters.
These parameters ensure that the estimated covariance matrix provides a lower and upper bound on the eigenvalues of the unknown $\bs{\Sigma}$.
This approach to defining the uncertainty sets of $\bs{\beta}$ and $\bs{\Sigma}$ in \eqref{def:B&E} is similar in spirit to \citep{delage2010_distr-rob-opt-mom-uncer}.
However, the description of $\mathcal{E}$ here also uses a matrix lower bound constraint on $\bs{\Sigma}$.
Another important motivation of introducing ambiguity set for the coefficient vector $\bs{\beta}^{ij}$
is on the linear utility model assumption. In principle, the utility can be a nonlinear function, and hence the usage of linear utility model can be biased. Imposing an ambiguity set on the linear coefficient vector can be viewed as a remedy 
for offsetting the possibly biased assumption of linearity.

\subsection{An Illustrative Example}
\label{sec:illus-exp}
We now provide a numerical example to illustrate that the choice of parameters 
in the specification of $\Puy$ and \eqref{def:B&E} may result in different decision recommendations. 
Consider a case that has three potential customer sites 
$S=\{1,2,3\}$ and three service center locations $F=\{1,2,3\}$.
The cost of opening a service center at each location is equal, and the service centers have unlimited capacity. The budget allows for opening
only one service center. Let the demand be $d_1=20$, $d_2=30$, and $d_3=25$.
Let $\bs{A}^{ij}=\textrm{diag}(a^{ij}_1,a^{ij}_2,a^{ij}_3)$ and 
$\whb{\Sigma}^{ij}=\textrm{diag}(\sigma^{ij}_1,\sigma^{ij}_2,\sigma^{ij}_3)$, 
i.e., the matrices $\bs{A}^{ij}$ and $\whb{\Sigma}^{ij}$ are diagonal for each $i,j\in\{1,2,3\}$. 
The model is given as follows:
\begin{displaymath}
\begin{aligned}
&\underset{\bs{y},\bs{x}}{\textrm{max}}\;\; \sum_{i,j\in\{1,2,3\}}\left( \underset{P\in\Puyij}{\textrm{min}}\;\E_P[u_{ij}(\bs{y})]\right)x_{ij} \\
&\textrm{ s.t. }\; y_1+y_2+y_3 \le 1,  \\
&\qquad  x_{1j}+x_{2j}+x_{3j} \le (d_1+d_2+d_3)y_j \qquad \forall j\in\{1,2,3\}, \\
&\qquad x_{i1} + x_{i2} +x_{i3} \le d_i \qquad\qquad\qquad\qquad \forall i\in\{1,2,3\}, \\
&\qquad x_{ij}\ge 0, \; y_j\in\{0,1\} \qquad\qquad\qquad\qquad \forall i,j\in\{1,2,3\}.
\end{aligned}
\end{displaymath}
In this case, there are only three possible decisions of the service center locations, which are \newline
$\bs{y}=[1,\;0,\;0],\;[0,\;1,\;0],\;[0,\;0,\;1]$.
We will see in Lemma~\ref{lem:a<E<b} that the optimal value $U^{ij}(\bs{y})=\underset{P\in\Puy}{\textrm{min}}\;\E_P[u_{ij}(\bs{y})]$
is given by the following form:
\begin{equation}\label{eqn:Uij-cond}
U^{ij}(\bs{y})=\textrm{max}\left\{
			(\hbeta^{ij})^{\top}\bs{y}-b^{ij}\|(\bs{A}^{ij})^{-1/2}\bs{y}\|,\; (\hbeta^{ij})^{\top}\bs{y}-\sqrt{\gamma^{ij}_2}\|(\whb{\Sigma}^{ij})^{1/2}\bs{y}\|
			\right\}.
\end{equation}
Suppose that the estimations of $\hat{\beta}^{ij}$ are given as follows:
\begin{displaymath}
\begin{aligned}
& \hat{\beta}^{11}=[8.5, 0.2, 0.4], \quad \hat{\beta}^{12}=[0.1, 8.0, 0.3], \quad \hat{\beta}^{13}=[0.2, 0.1, 7.3], \\
& \hat{\beta}^{21}=[8.2, 0.0, 0.2], \quad  \hat{\beta}^{22}=[0.1, 8.2, 0.3], \quad  \hat{\beta}^{23}=[0.2, 0.0, 7.4], \\
& \hat{\beta}^{31}=[8.3, 0.1, 0.2], \quad  \hat{\beta}^{32}=[0.0, 8.1, 0.1], \quad \hat{\beta}^{33}=[0.1, 0.0, 7.5].
\end{aligned}
\end{displaymath}

In the base case we assume that all parameter estimates are exact, i.e., there is no ambiguity. 
\begin{displaymath}
\begin{aligned}
&a^{ij}_k = \sigma^{ij}_k = \gamma^{ij}_2 = 0 \quad \forall i,j,k\in\{1,2,3\}, \quad b^{i1} = 0, \quad b^{i2}=0, \quad b^{i3}=0 \quad \forall i\in \{1,2,3\}.
\end{aligned}
\end{displaymath}
Now consider two different parameter settings for  $\bs{A}^{ij}$ and $\widehat{\bs{\Sigma}}^{ij}$.\newline
\noindent Parameter Estimation 1:
\begin{displaymath}
\begin{aligned}
&a^{ij}_k = \sigma^{ij}_k = \gamma^{ij}_2 = 2.0 \quad \forall i,j,k\in\{1,2,3\}, \quad b^{i1} = 1.41, \quad b^{i2}=1.27, \quad b^{i3}=2.69 \quad \forall i\in \{1,2,3\}.
\end{aligned}
\end{displaymath}
\noindent Parameter Estimation 2:
\begin{displaymath}
\begin{aligned}
&a^{ij}_k = \sigma^{ij}_k = \gamma^{ij}_2 = 2.0 \quad \forall i,j,k\in\{1,2,3\}, \quad b^{i1} = 1.41, \quad b^{i2}=0.99, \quad b^{i3}=2.55 \quad \forall i\in \{1,2,3\}.
\end{aligned}
\end{displaymath} 
We can verify that in the base case $\bs{y}^*=[1,\;0,\;0]$ is the optimal solution with the optimal value 582. This is also the solution, with the optimal value 548, under the parameter Estimation 1. However, $\bs{y}^*=[0,\;1,\;0]$ is the optimal solution with the optimal value 556 under the parameter Estimation 2. Note that in comparison to Estimation 1, the level of ambiguity at locations 2 and 3 is smaller (the $b^{i2}$ and $b^{i3}$ values are smaller) 
in Estimation 2. This reduced ambiguity results in a different service center location decision.
We arrive at different decisions with varying levels of ambiguity.

\subsection{Potential simplification}
A shortcoming of the $\bs{\beta}^{ij}$ vector estimation model \eqref{eqn:linear-reg}
is that for each $(i,j)$ pair, some kind of questionnaire is required to conduct in order to
collect samples needed to fit the linear regression model. Practically, this way of sample
collection may not be scalable. One way to make the sample collection and coefficient estimation more efficient is
to interpret $\bs{\beta}^{ij}$ as the multiplication of a $(i,j)$ dependent matrix with
a $(i,j)$ independent vector, i.e., $\bs{\beta}^{ij}=M^{ij}\bs{\beta}$. Under this setting,
the estimation of the vector $\bs{\beta}$ can be based on samples collected that are 
not restricted to a specific $(i,j)$ pair if it is assumed that customers from all residential sites 
are homogeneous and the attractiveness of a facility is mainly determined by some standard
factors such as travel distance, scale and location environment, etc. The $M^{ij}$ can be a 
pre-determined 0-1 matrix that categorizes each candidate location in the neighborhood of 
the residential sit $i$. With this simplification, the vector $\bs{\beta}$ will be used for modeling 
the utility gain for all $(i,j)$ pairs and the ambiguity set only needs to be constructed for
a single vector $\bs{\beta}$.

Another way of simplification is to assume the utility coefficient vector only depend on the customer 
site $i\in S$, i.e., the set of utility coefficient vectors are $\{\bs{\beta}^i:\;i\in S\}$. The heterogeneity
of customer groups can be modeled in this setting, and it can be more practical to collect data to fit the utility
models than for the case of having $(i,j)$-dependent coefficient vectors.

\section{Mixed 0-1 Conic Reformulation of \eqref{opt:RFL}}
\label{sec:mixed-bin-conic}
We first give an analytical solution of the inner problem of \eqref{opt:RFL-II}.
We show that the analytical solution can be written as the maximum of two second-order-cone functions. 
This analytical solution is used to reformulate \eqref{opt:RFL-II} as a mixed 0-1 second-order-cone program.
\subsection{Reformulation of \eqref{opt:RFL} Using the Moment Based Ambiguity Set}
\label{sec:RFL-reform}
We first reformulate the inner problem in \eqref{opt:RFL-II}:
\begin{equation}\label{opt:inner-ij}
\textrm{min}_{P\in\Puy}\;\E_P[u(\bs{y})],
\end{equation}
where $\Puy$ is defined in \eqref{def:B&E}. We have omitted the subscripts $i,j$ for simplicity. The following proposition provides an explicit decision dependent description of the ambiguity set $\Puy$.

\begin{proposition}\label{prop:Puy-reform}
Let $\Puy$ be defined as in \eqref{eqn:amb-P} and \eqref{def:B&E}. Then $\Puy$ can be reformulated as:
\begin{equation}\label{eqn:Puy-reform}
\Puy=\Set*{P\in(\Xi,\mtc{F})}{\begin{aligned}
                                & \hbeta^{\top}\bs{y}-b\|\bs{A}^{-1/2}\bs{y}\| \le\E_P[u]\le\hbeta^{\top}\bs{y}+b\|\bs{A}^{-1/2}\bs{y}\| \\
                                & \gamma_1\bs{y}^{\top}\whb{\Sigma}\bs{y}\le \E_{P}[(u-\hbeta^{\top}\bs{y})^2]\le\gamma_2\bs{y}^{\top}\whb{\Sigma}\bs{y} 
                               \end{aligned}}.
\end{equation}
\end{proposition}
\begin{proof}
We first show the reformulation of $\E_P[u]$. Inequality \eqref{def:B&E} that defines $\mtc{B}$ is equivalent to $\|\bs{A}^{1/2}(\bs{\beta}-\hbeta)\|\le b$.
Then the following upper bound on $\bs{\beta}^{\top}\bs{y}$ holds for any decision vector $\bs{y}$.
\begin{displaymath}
\begin{aligned}
&\bs{\beta}^{\top}\bs{y}=\hbeta^{\top}\bs{y}+(\bs{\beta}-\hbeta)^{\top}\bs{y}
=\hbeta^{\top}\bs{y}+(\bs{\beta}-\hbeta)^{\top}\bs{A}^{1/2}\bs{A}^{-1/2}\bs{y} \\
&\le\hbeta^{\top}\bs{y}+\|\bs{A}^{1/2}(\bs{\beta}-\hbeta)\|\|\bs{A}^{-1/2}\bs{y}\|
\le\hbeta^{\top}\bs{y}+b\|\bs{A}^{-1/2}\bs{y}\|,
\end{aligned}
\end{displaymath}
where we use the Cauchy-Schwarz inequality and the definition of the ellipsoid
in the first and the second inequalities of the above expression, respectively. 
Similarly, we can show that $\bs{\beta}^{\top}\bs{y}\ge\hbeta^{\top}\bs{y}-b\|\bs{A}^{-1/2}\bs{y}\|$.
Furthermore, the above upper and lower bounds on $\bs{\beta}^{\top}\bs{y}$ are atttainable based on
the conditions for the equality to hold in the Cauchy-Schwarz inequality. Therefore, we have
\begin{equation}\label{eqn:minmax-beta*y}
\underset{\bs{\beta}\in\mtc{B}}{\textrm{max}}\;\bs{\beta}^{\top}\bs{y} = \hbeta^{\top}\bs{y}+b\|\bs{A}^{-1/2}\bs{y}\|,
\qquad
\underset{\bs{\beta}\in\mtc{B}}{\textrm{min}}\;\bs{\beta}^{\top}\bs{y} = \hbeta^{\top}\bs{y}-b\|\bs{A}^{-1/2}\bs{y}\|.
\end{equation}
Based on the definition of $\mathcal{E}$ in \eqref{def:B&E}, we have 
\begin{equation}\label{eqn:minmax-y*sigma*y}
\underset{\bs{\Sigma}\in\mathcal{E}}{\textrm{min}}\;\bs{y}^{\top}\bs{\Sigma}\bs{y}=\gamma_1\bs{y}^{\top}\whb{\Sigma}\bs{y},
\qquad \underset{\bs{\Sigma}\in\mathcal{E}}{\textrm{max}}\;\bs{y}^{\top}\bs{\Sigma}\bs{y}=\gamma_2\bs{y}^{\top}\whb{\Sigma}\bs{y}.
\end{equation} 
Substituting \eqref{eqn:minmax-beta*y} and \eqref{eqn:minmax-y*sigma*y} into \eqref{eqn:amb-P},
we obtain \eqref{eqn:Puy-reform}.
\end{proof}
The following lemma allows us to solve \eqref{opt:inner-ij} analytically.
\begin{lemma}\label{lem:a<E<b}
Let $\Xi=[a,b]$ be a closed interval in $\mathbb{R}$ and $\mathcal{F}$ be the Borel $\sigma$-algebra on $\Xi$. 
Let $(\Xi,\mathcal{F})$ denote the set of probability measures defined on $\Xi$ with the Borel $\sigma$-algebra.
Consider the following optimization problem:
\begin{equation}\label{opt:minEp[u]}
\emph{min}_{P\in(\Xi,\mathcal{F})}\;\E_P[u]\qquad \emph{s.t.}\quad c_1\le \E_P[u]\le c_2, \quad d_1\le\E_P[(u-\mu)^2]\le d_2,
\end{equation}
where $u$ is a random variable and $a_1,a_2,d_1,d_2,\mu\in\mathbb{R}$ 
satisfying $a\le c_1\le\mu\le c_2\le b$, $0\le d_1\le d_2$ and $(\mu-c_1)^2\ge d_1$,
which guarantees that problem \eqref{opt:minEp[u]} is feasible.
Let $V^*$ be the optimal value of the problem. Then we have $V^*=\emph{max}\{c_1,\;\mu-\sqrt{d_2}\}$.
\end{lemma}
\begin{proof}
We consider the solution of \eqref{opt:minEp[u]} in two cases: $(\mu-c_1)^2\le d_2$ and $(\mu-c_1)^2> d_2$.
If $(\mu-c_1)^2\le d_2$, one can construct a probability measure $P^*$ such that $P^*(\{c_1\})=1$.
The measure $P^*$ is feasible since $\E_{P^*}[u]=c_1$ and $\E_P[(u-\mu)^2]=(\mu-c_1)^2\in [d_1,d_2]$ satisfying the constraints. It is also optimal since $\E_{P^*}[u]=c_1$ by construction and hence satisfying the lower bound constraint $\E_P[u]\ge c_1$ at equality.
In this case, we have $V^*=c_1$. The conditions $c_1\le \mu$ and $(\mu-c_1)^2\le d_2$ further imply that 
$c_1>\mu-\sqrt{d_2}$. Therefore, in this case the expression $V^*=\textrm{max}\{c_1,\;\mu-\sqrt{d_2}\}$ holds.

Now consider the case that $(\mu-c_1)^2> d_2$. Due to the constraint on the second moment of $u$, we have
$\int_{[a,b]}(u-\mu)^2dP(u)\le d_2$. The Cauchy-Schwarz inequality gives 
\begin{displaymath}
\int_{[a,b]}(u-\mu)^2dP(u)\ge \left(\int_{[a,b]}(u-\mu)dP(u)\right)^2=\left(\E_P[u]-\mu\right)^2,
\end{displaymath}  
and hence $d_2\ge\left(\E_P[u]-\mu\right)^2$, which implies that $\sqrt{d_2}\ge \mu-\E_P[u]$.
Therefore, we have $V^*\ge\mu-\sqrt{d_2}$. It remains to show that the lower bound $\mu-\sqrt{d_2}$
is attainable. We now construct an optimal probability measure $P^*$ such that $P^*(\{\mu-\sqrt{d_2}\})=1$.
To verify that $P^*$ is feasible, we note that $\E_{P^*}[u]=\mu-\sqrt{d_2}$ which is in the interval $[c_1,c_2]$
due to the conditions $\mu\le c_2$ and $(\mu-c_1)^2>d_2$. Furthermore, we have $\E_{P^*}[(u-\mu)]=d_2$.
Combining the above two cases, we get $V^*=\textrm{max}\{c_1,\;\mu-\sqrt{d_2}\}$ which concludes the proof.
\end{proof}
\begin{remark}
According to Lemma~\ref{lem:a<E<b} the optimal value of \eqref{opt:minEp[u]} does not depend on the value of constants $c_2$ and $d_1$ in the constraints. We now provide an interpretation
of the optimal value. In the case of $(\mu-c_1)^2\le d_2$, the maximum deviation allows the mean value 
$\mathbb{E}_P[u]$ to reach the lower bound. In the case of $(\mu-c_1)^2> d_2$, the deviation is not large enough.
Consequently, the lower bound of the mean value is not attainable, and the optimal value depends on the maximum deviation determined by $d_2$. 
Note that allowing a matrix lower bound in the definition of $\mathcal{E}$ is different from the setting (1b) 
of the moment-based ambiguity set in \citep{delage2010_distr-rob-opt-mom-uncer}. If a matrix lower bound 
is imposed in (1b) of \citep{delage2010_distr-rob-opt-mom-uncer}, the distributionally-robust optimization model
in \citep{delage2010_distr-rob-opt-mom-uncer} can not be reformulated into a convex optimization problem.
However, an analytical specification of the optimal value of \eqref{opt:minEp[u]} is possible because the ambiguity set is defined for
the univariate utility and simplification is possible in this case.
By applying Lemma~\ref{lem:a<E<b} to \eqref{eqn:Puy-reform}, it gives an analytical expression for the 
optimal value of \eqref{opt:inner-ij}. The result is given in Corollary~\ref{cor:anal-value}.
\end{remark}
\begin{corollary}\label{cor:anal-value}
The optimal value of \eqref{opt:inner-ij} is given as follows:
\begin{equation}\label{eqn:sol-inner}
\emph{min}_{P\in\Puyij}\;\E_P[u_{ij}(\bs{y})]=\emph{max}\left\{
			(\hbeta^{ij})^{\top}\bs{y}-b^{ij}\|(\bs{A}^{ij})^{-1/2}\bs{y}\|,\; (\hbeta^{ij})^{\top}\bs{y}-\sqrt{\gamma^{ij}_2}\|(\whb{\Sigma}^{ij})^{1/2}\bs{y}\|
			\right\},
\end{equation}
for any $i\in S$ and $j\in F$.
\end{corollary}

When substituting the optimal value \eqref{eqn:sol-inner} into \eqref{opt:RFL-II},
we get a nonlinear term written as $x_{ij}\left(\textrm{min}_{P\in\Puyij}\;\E_P[u_{ij}(\bs{y})]\right)$.
This nonlinear term involves bilinear product terms $x_{ij}\bs{y}$. 
A reformulation of \eqref{opt:RFL-II} based on linearizing these bilinear product terms is given in the following proposition.
\begin{proposition}\label{prop:inner-reform}
Let the ambiguity set $\Puyij$ be defined as in \eqref{eqn:amb-P}--\eqref{def:B&E} for all $i\in S$ and $j\in F$.
The recourse problem \eqref{opt:RFL-II} can be reformulated as 
\begin{equation}\label{opt:RFL-II-0}
\begin{aligned}
&\emph{max}\;\; \sum_{i\in S}\sum_{j\in F}\Uij \\
&\;\emph{ s.t.}\;\; \Uij\le\emph{max}\{(\hbeta^{ij})^{\top}\bs{v}^{ij}-b^{ij}\|(\bs{A}^{ij})^{-1/2}\bs{v}^{ij}\|,
                                  \;  (\hbeta^{ij})^{\top}\bs{v}^{ij}-\sqrt{\gamma^{ij}_2}\|(\whb{\Sigma}^{ij})^{1/2}\bs{v}^{ij}\|\}  \\
&\qquad\qquad  \forall i\in S,\; \forall j\in F, \\
&\qquad v^{ij}_k\le R^{ij}y_k, \quad v^{ij}_k\le x_{ij},  
	\quad v^{ij}_k\ge x_{ij}-R^{ij}(1-y_k) \quad \forall i\in S,\; \forall j\in F,\; \forall k\in F, \\
&\qquad \bs{x}\in X(\bs{y}),\;U^{ij}\ge 0,\; \bs{v}^{ij}\in\mathbb{R}^{|F|}_+ \hspace{4.5cm}\forall i\in S,\;\forall j\in F,
\end{aligned}
\tag{RSP-0}
\end{equation}
where $\Rij$ is a constant satisfying $R^{ij}=\emph{min}\{D_i,C_j\}$.  
\end{proposition}
\begin{proof}
By Lemma~\ref{lem:a<E<b}, the optimal value of \eqref{opt:inner-ij} can be written as
\begin{equation}
\textrm{max}\{(\hbeta^{ij})^{\top}\bs{y}x_{ij}-b^{ij}\|(\bs{A}^{ij})^{-1/2}\bs{y}x_{ij}\|,
                                  \;\;  (\hbeta^{ij})^{\top}\bs{y}x_{ij}-\sqrt{\gamma^{ij}_2}\|(\whb{\Sigma}^{ij})^{1/2}\bs{y}x_{ij}\| \}.
\end{equation}
We need to verify that the second to the fourth constraints in \eqref{opt:RFL-II-0} ensure that $v^{ ij}_k=x_{ij}y_k$
for all $i\in S$, $j\in F$ and $k\in F$. If $y_k=0$, the second constraint implies that $v^{ ij}_k\le 0$. Combining it with the non-negative
constraint on $v^{ ij}_k$ implies that $v^{ ij}_k=0$. 
If $y_k=1$, the third and fourth constraints imply that $v^{ ij}_k=x_{ij}$.
Therefore, the recourse problem \eqref{opt:RFL-II} can be reformulated as \eqref{opt:RFL-II-0}.
\end{proof}
It is worth to remark that a constraint of the type $x^2_3\le x_1x_2$ which involves a bilinear
term can be reformulated as a second-order-cone constraint $\norm{(2x_3, x_1-x_2)}\le x_1+x_2$. 
This technique plays an essential role in reformulating multinomial logit choice model based constraints
into compact mixed-integer second-order-cone constraints 
\citep{sen2018-conic-IP-assort-multi-logit,lin2020-opt-locker-loc-logit-choice-model}. 
Unfortunately this technique is not applicable to simplify our reformulation.

Another approach for strengthening 
the mixed-binary linear representation of the bilinear terms 
in \eqref{opt:RFL-II-0} is to use McCormick envelopes.
Ignoring indices $i,j$ for now, we have define auxiliary variables $v_k$ to represent the bilinear scalar $y_kx$
and perform the linearization $v_k\le Ry_k$, $v_k\le x$, $v_k\ge x-R(1-y_k)$ which
provides a necessary and sufficient condition to ensure that $v_k=y_kx$ holds. 
The linearization can be potentially strengthened
by adding the following McCormick-envelope valid inequalities \citep{mccormick1976} to the reformulation:
$v_k\ge x^Ly_k+xy^L-x^Ly^L$, $v_k\ge x^Uy_k + xy^U - x^Uy^U$,
$v_k\le x^Uy_k + xy^L - x^Uy^L$ and $v_k\le xy^U + x^Ly_k - x^Ly^U$,
where $x^L$, $x^U$ (resp. $y^L$, $y^U$) are lower and upper bounds for $x$ (resp. $y_k$).
Using $x^L=0$, $y^L=0$, $x^U=R$, $y^U=1$, the above inequalities become
$v_k\ge Ry_k+x-R$, $v_k\le Ry_k$ and $v_k\le x$, which is exactly the same set of inequalities
as from the linearization. This means in our case, the McCormick envelopes trivially reduce to the
linearization inequalities.

The most challenging part in \eqref{opt:RFL-II-0} is the first constraint. 
Note that the two functions inside the `max' are both concave, and hence the first constraint is non-convex.
We can reformulate this non-convex constraint into convex constraints by introducing binary variables in the model.
This reformulation is given in Section~\ref{sec:reform-lift}. 
Here we highlight an important special case of problem
setting under which \eqref{opt:RFL-II-0} can be simplified significantly. Consider the special case for which
the condition $(\bs{A}^{ij})^{-1/2}=(\whb{\Sigma}^{ij})^{1/2}$ is satisfied for all $i$ and $j$.
In this case, the first non-convex constraint of \eqref{opt:RFL-II-0} can be equivalently formulated as
$$\Uij\le(\hbeta^{ij})^{\top}\bs{v}^{ij}-\min\{b^{ij}, \sqrt{\gamma^{ij}_2}\}\|(\whb{\Sigma}^{ij})^{1/2}\bs{v}^{ij}\|,$$
which is clearly a SOC constraint. 

\subsection{Reformulation using Convexification in a Lifted Space}
\label{sec:reform-lift}
In this reformulation of \eqref{opt:RFL-II-0}, we lift the feasible set of the variables $\{U,\bs{v}\}$ (omitting indices $i,j$)
into a higher dimensional space represented by variables $\{U,U^1,U^2,\bs{v}^1,\bs{v}^2, \bs{s}\}$, 
where $\{U^1,\bs{v}^1\}$, $\{U^2,\bs{v}^2\}$ are additional variables
introduced to represent the max constraint in \eqref{opt:RFL-II-0},
and $\bs{s}$ are selection variables that are used to determine which
subset of variables (either $\{U^1,\bs{v}^1\}$ or $\{U^2,\bs{v}^2\}$) are active.
This reformulation does not use a big-M constant 
to handle the ``max'' involved in the first constraint with the help of a lifting technique.
A reformulation is also possible using big-M constants without lifting, 
but it is omitted here because its performance was not superior to the one given here.
Note that we still need to use big-M coefficients $R^{ij}$'s to reformulate
the bilinear terms $x_{ij}y_k$, but it is a separate issue.
\begin{theorem}\label{thm:RFL-II-2}
The recourse problem \eqref{opt:RFL-II} can be reformulated as the following mixed 0-1 second-order-cone programming problem:
\begin{equation}\label{opt:RFL-II-2}
\begin{aligned}
& \emph{max}\;\; \sum_{i\in S}\sum_{j\in F}\Uij \\
& \emph{ s.t. }\; \Uijo \le (\hbeta^{ij})^{\top}\bvijo-b^{ij}\|(\bs{A}^{ij})^{-1/2}\bvijo\|,\quad \Uijt \le (\hbeta^{ij})^{\top}\bvijt-\sqrt{\gamma^{ij}_2}\|(\whb{\Sigma}^{ij})^{1/2}\bvijt\|, \\
& \qquad\qquad \forall i\in S,\;\forall j\in F, \\
&\qquad \Uij = \Uijo+\Uijt, \quad \bvij = \bvijo+\bvijt, \quad v^{ij1}_k\le R^{ij}\sij,
   \quad v^{ij2}_k\le R^{ij}(1-\sij)    \\
&  \qquad\qquad \forall i\in S,\;\forall j\in F,\;\forall k\in F, \\
&\qquad v^{ij}_k\le R^{ij}y_k,\quad v^{ij}_k\le x_{ij},\quad v^{ij}_k\ge x_{ij}-R^{ij}(1-y_k),
    \quad \forall i\in S,\;\forall j\in F,\;\forall k\in F, \\
&\qquad \bs{x}\in X(\bs{y}),\; s^{ij}\in\{0,1\},\;  U^{ij}, \Uijo, \Uijt \ge 0, 
 \bs{v}^{ij}, \bvijo, \bvijt \in\mathbb{R}^{|F|}_+ \quad \forall i\in S,\;\forall j\in F. 
\end{aligned}
\tag{RSP-1}
\end{equation}
\end{theorem}
\begin{proof}
It is easy to see that if $s^{ij}=0$, $v^{ij1}_k=0$, $U^{ij1}=0$, $\bs{v}^{ij}=\bs{v}^{ij2}$,
and $U^{ij}=U^{ij2}=(\hbeta^{ij})^{\top}\bvijt-\sqrt{\gamma^{ij}_2}\|(\whb{\Sigma}^{ij})^{1/2}\bvijt\|$.
If $s^{ij}=1$, $v^{ij2}_k=0$, $U^{ij2}=0$, $\bs{v}^{ij}=\bs{v}^{ij1}$,
and $U^{ij}=U^{ij1}=(\hbeta^{ij})^{\top}\bvijo-b\|(\bs{A}^{ij})^{-1/2}\bvijo\|$.
Therefore, the recourse problem \eqref{opt:RFL-II} can be reformulated into \eqref{opt:RFL-II-2}.
\end{proof}

\begin{corollary}
Let the ambiguity set $\Puyij$ be defined as \eqref{eqn:amb-P} and \eqref{def:B&E},  
the \eqref{opt:RFL} is reformulated as:
\begin{equation}\label{opt:RFL-determ}
\begin{aligned}
& \emph{max}\;\; \bs{h}^{\top}\bs{y}+\sum_{i\in S}\sum_{j\in F}U^{ij} \\
& \emph{ s.t. }\; U^{ij1} \le (\hbeta^{ij})^{\top}\bs{v}^{ij1}-b^{ij}\|(\bs{A}^{ij})^{-1/2}\bs{v}^{ij1}\|,\quad
U^{ij2} \le (\hbeta^{ij})^{\top}\bs{v}^{ij2}-\sqrt{\gamma^{ij}_2}\|(\whb{\Sigma}^{ij})^{1/2}\bs{v}^{ij2}\| \\
&\qquad\qquad \forall i\in S,\;\forall j\in F, \\
&\qquad U^{ij} = U^{ij1}+U^{ij2}, \quad \bs{v}^{ij} = \bs{v}^{ij1}+\bs{v}^{ij2},\quad v^{ij1}_k\le R^{ij}s^{ij},  
   \quad v^{ij2}_k\le R^{ij}(1-s^{ij}) \\
&\qquad\qquad \forall i\in S,\;\forall j\in F,\;\forall k\in F, \\
&\qquad v^{ij}_k\le R^{ij}y_k,\quad v^{ij}_k\le x_{ij},\quad v^{ij}_k\ge x_{ij}-R^{ij}(1-y_k) \quad \forall i\in S,\;\forall j\in F,\;\forall k\in F, \\
&\qquad \bs{x}\in X(\bs{y}),\quad \sum_{j\in F}b_jy_j\le B, \\
&\qquad  U^{ij}, U^{ij1}, U^{ij2} \ge 0, 
\bs{v}^{ij}, \bs{v}^{ij1}, \bs{v}^{ij2} \in\mathbb{R}^{|F|}_+ \;\; \bs{y},\bs{s}\in\{0,1\}^{|F|}\;\; \forall i\in S,\;\forall j\in F,
\end{aligned}
\tag{RFL-MISOCP}
\end{equation}
where $X(\bs{y})$ is defined in \eqref{eqn:X-set}.
\end{corollary}

\section{Generating a Convex Hull of the Max Substructure}
\label{sec:gen-conv-hull}
A challenge in solving \eqref{opt:RFL} comes from the max inequality in \eqref{opt:RFL-II-0}. 
This inequality is re-written as 
\begin{equation}\label{eqn:nonconv-constr}
\Uij\le\textrm{max}\{f^{ ij}(\bs{v}^{ij}),\;g^{ij}(\bs{v}^{ij})\},
\end{equation}
where $f^{ij}(\bs{v})=(\hbeta^{ij})^{\top}\bs{v}-b^{ij}\|(\bs{A}^{ij})^{-1/2}\bs{v}\|$, 
and $g^{ij}(\bs{v})=(\hbeta^{ij})^{\top}\bs{v}-\sqrt{\gamma^{ij}_2}\|(\whb{\Sigma}^{ij})^{1/2}\bs{v}\|$.
Note that the second-order-conic functions $f^{ij}(\bs{v})$ and $g^{ij}(\bs{v})$ are concave, and therefore, the maximal function on the right side of \eqref{eqn:nonconv-constr} is not concave and the constraint 
\eqref{eqn:nonconv-constr} gives a non-convex feasible set. 
The reformulation \eqref{opt:RFL-II-2} introduces 
extra binary (continuous) variables and constraints to 
reformulate this non-convex constraint based region as mixed-binary conic constraints, 
whose relaxation is a convex set. 
This reformulation introduces a large amount of auxiliary (binary and continuous) variables 
and constraints, which could be less efficient for medium and large problem instances. 
An alternative approach is to avoid the challenge by considering a relaxed formulation 
of \eqref{opt:RFL-II-0} as follows:
\begin{equation}\label{opt:RSP-relax}
\begin{aligned}
&\text{max}\;\; \sum_{i\in S}\sum_{j\in F}\Uij & \\
&\;\text{ s.t. } \Uij\le \psi^{ij}(\bs{v}^{ij}) &  \forall i\in S,\; \forall j\in F, \\
&\qquad\; v^{ij}_k\le R^{ij}y_k, \quad v^{ij}_k\le x_{ij},  
	\quad v^{ij}_k\ge x_{ij}-R^{ij}(1-y_k) & \forall i\in S,\; \forall j\in F,\; \forall k\in F, \\
&\qquad\; \bs{x}\in X(\bs{y}),\;U^{ij}\ge 0,\; \bs{v}^{ij}\in\mathbb{R}^{|F|}_+ & \forall i\in S,\;\forall j\in F, 
\end{aligned}
\tag{RSP-relax}
\end{equation}
where $\psi^{ij}$ is the function such that the subgraph of $\psi^{ij}$ is exactly the convex hull
of the following set
\begin{equation}
\chi^{ij} = \Set*{(U,\bs{v})}{ U\le \textrm{max}\{f^{ij}(\bs{v}),g^{ij}(\bs{v})\}}.
\end{equation}
Substituting the recourse problem relaxation \eqref{opt:RSP-relax} into \eqref{opt:RFL},
we end up with the following relaxation of \eqref{opt:RFL}:
\begin{equation}\label{opt:RFL-relax}
\begin{aligned}
&\text{max}\;\; \bs{h}^\top\bs{y}+\sum_{i\in S}\sum_{j\in F}\Uij & \\
&\;\text{ s.t. } \Uij\le \psi^{ij}(\bs{v}^{ij}) &  \forall i\in S,\; \forall j\in F, \\
&\qquad\; v^{ij}_k\le R^{ij}y_k, \quad v^{ij}_k\le x_{ij},  
	\quad v^{ij}_k\ge x_{ij}-R^{ij}(1-y_k) & \forall i\in S,\; \forall j\in F,\; \forall k\in F, \\
&\qquad\; \bs{x}\in X(\bs{y}),\;U^{ij}\ge 0,\; \bs{v}^{ij}\in\mathbb{R}^{|F|}_+ & \forall i\in S,\;\forall j\in F, \\
&\qquad\; \bs{y}\in\{0,1\}^{|F|}.
\end{aligned}
\tag{RFL-relax}
\end{equation}
We develop a cutting-plane algorithm to solve \eqref{opt:RFL-relax} iteratively. 
The core part of the algorithm is to generate tangent cutting planes to approximate
the function $\chi^{ij}$ with a finite number of hyper planes as needed. 
Note that in general, solving \eqref{opt:RFL-relax} to optimality may not 
generate an optimal solution of \eqref{opt:RFL} due to the relaxation, 
but it should yield a high-quality solution.

\subsection{Strengthen Formulations using Tangent Planes}
\label{sec:conv-1}
We omit the indices $\omega,i,j$ to simplify the notations in the following discussion. 
The convex hull $\conv(\chi)$ can be re-written as $\conv(\chi)=\conv(\chi_1\cup\chi_2)$,
where $\chi_1$ and $\chi_2$ are the following convex sets:
\begin{equation}
\chi_1=\Set*{(U,\bs{v})}{0\le U\le f(\bs{v}),\;\bs{v}\ge 0}, \qquad \chi_2=\Set*{(U,\bs{v})}{0\le U\le g(\bs{v}),\;\bs{v}\ge 0}.
\end{equation}
To describe $\conv(\chi)$, it suffices to provide all tangent inequalities of $\conv(\chi)$.
\begin{definition}\label{def:facet}
Let $\mathcal{C}$ be a convex set in $\mathbb{R}^n$. 
A linear inequality $\bs{a}^{\top}\bs{x}\le b$ is a \emph{tangent inequality} of $\mathcal{C}$ if 
the inequality is valid for all points in $\mathcal{C}$ and the intersection set
$\Set*{\bs{x}\in\mathbb{R}^n}{\bs{a}^{\top}\bs{x}=b}\cap\mathcal{C}$ is non-empty.
For a tangent inequality $h$ of $\mathcal{C}$ represented by the inequality $\bs{a}^{\top}\bs{x}\le b$, the intersection set 
$\Set*{\bs{x}\in\mathbb{R}^n}{\bs{a}^{\top}\bs{x}=b}\cap\mathcal{C}$ is defined as the
\emph{tangent points} of $h$, and it is denoted by $\mathcal{T}(h,\mathcal{C})$.
\end{definition}

\begin{figure}
\centering
\includegraphics[scale=0.5, trim=1.5cm 15cm 0.5cm 4cm]{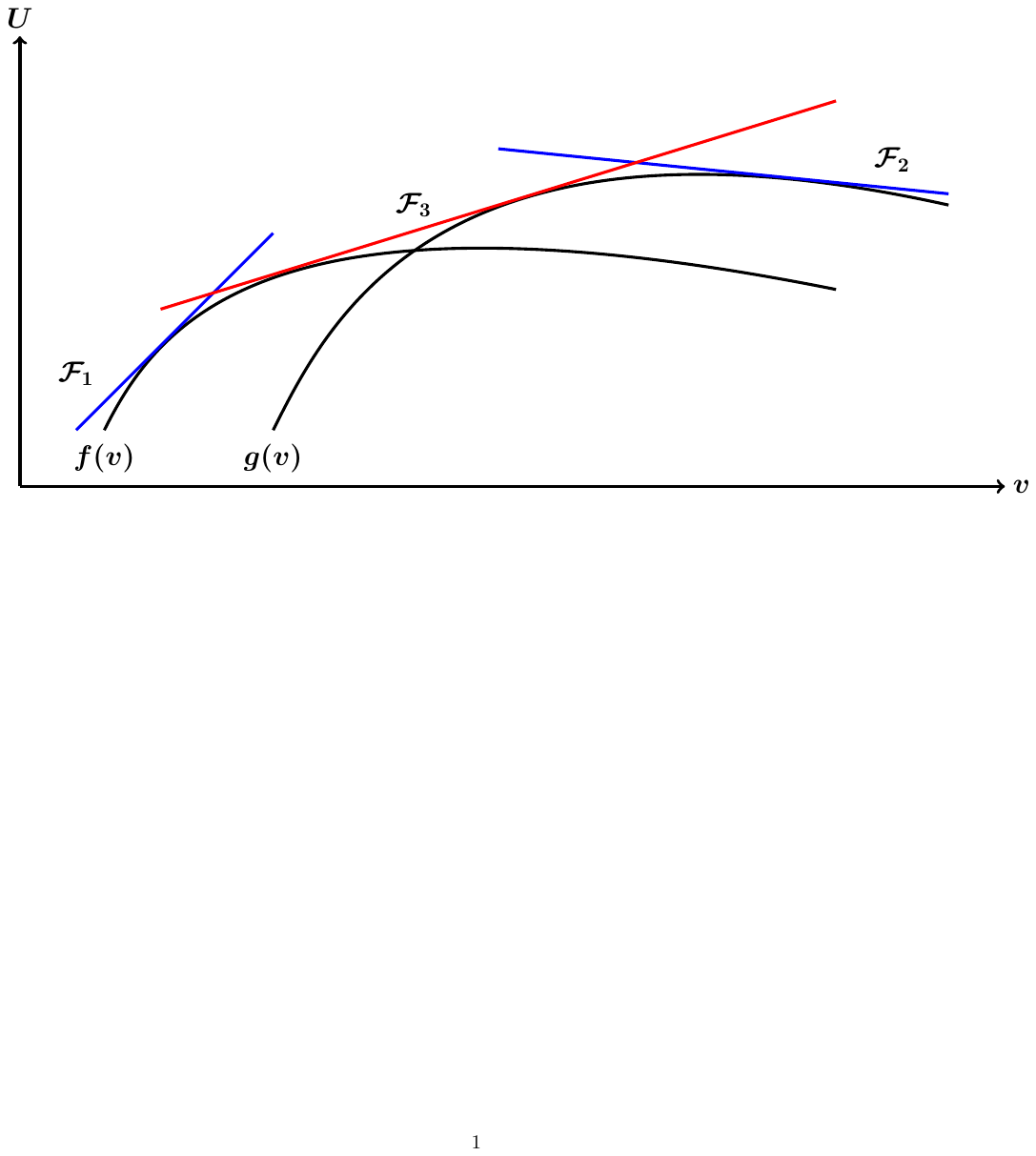}
\caption{Illustration of tangent inequalities in $\mathcal{F}_1$,  $\mathcal{F}_2$ and  $\mathcal{F}_3$. 
The tangent inequalities in $\mathcal{F}_1\cup\mathcal{F}_2$ are generated as G-cuts \eqref{eqn:ineq-grad1}-\eqref{eqn:ineq-grad2},
and tangent inequalities in $\mathcal{F}_3$ are generated as P-cuts \eqref{eqn:P-cut}.}\label{fig:illu-cuts}
\end{figure}

Note that the functions $f(\bs{v})$ and $g(\bs{v})$ in our case are differentiable everywhere except at the origin $(0,\bs{0})$.
At a point $\bs{v}\in\mathbb{R}^{|F|}_+\setminus\{(0,\bs{0})\}$, 
a tangent inequality of $\conv(\chi)$ corresponds to a tangent plane of $\conv(\chi)$. 
Let $\mF^0$ be the set of tangent inequalities of $\conv(\chi)$, and let $\mF$ be a subset of $\mF^0$ defined as follows:
\begin{equation}\label{def:F}
\mF=\Set*{p\in\mF^0}{(0,\bs{0})\notin\mathcal{T}(p,\conv(\chi))}.
\end{equation} 
We focus on investigating $\mF$ instead of $\mF^0$ to avoid dealing 
with the non-differentiable point $(0,\bs{0})$ in the discussion.
The point $(0,\bs{0})$ is handled in Theorem~\ref{thm:conv-hull-tang-plane}.
The tangent inequalities in $\mF$ can be partitioned into the following three disjointed subsets:
\begin{enumerate}
	\item $\mathcal{F}_1$: Tangent inequalities corresponding to hyperplanes that are only tangent to $\chi_1$;
	\item $\mathcal{F}_2$: Tangent inequalities corresponding to hyperplanes that are only tangent to $\chi_2$;
	\item $\mathcal{F}_3$: Tangent inequalities corresponding to hyperplanes that are tangent to both $\chi_1$ and $\chi_2$.
\end{enumerate}
The illustration of the tangent inequalities of $\conv(\chi_1\cup\chi_2)$ is given in Figure~\ref{fig:illu-cuts}.
\begin{proposition}\label{prop:grad-cut}
For any point $\bs{v}^{\prime}\in\mathbb{R}^{|F|}_+\setminus\{\bs{0}\}$, 
consider the following two convex optimization problems:
\begin{align}
&\psi_1(\bs{v}^{\prime})=\underset{\bs{v}}{\emph{min}}\;\nabla f(\bs{v}^{\prime})^{\top}(\bs{v}-\bs{v}^{\prime})+f(\bs{v}^{\prime})-g(\bs{v}), \label{opt:f-g} \\
&\psi_2(\bs{v}^{\prime})=\underset{\bs{v}}{\emph{min}}\;\nabla g(\bs{v}^{\prime})^{\top}(\bs{v}-\bs{v}^{\prime})+g(\bs{v}^{\prime})-f(\bs{v}). \label{opt:g-f}
\end{align}
The subsets $\mF_1$, $\mF_2$ and $\mF_3$ can be represented as follows:
\begin{equation}\label{eqn:F123}
\begin{aligned}
&\mF_1=\Set*{U\le \nabla f(\bs{v}^{\prime})^{\top}(\bs{v}-\bs{v}^{\prime})+f(\bs{v}^{\prime})}{\psi_1(\bs{v}^{\prime})>0,\; \bs{v}^{\prime}\in\mathbb{R}^{|F|}_+\setminus\{\bs{0}\}},\\
&\mF_2=\Set*{U\le \nabla g(\bs{v}^{\prime})^{\top}(\bs{v}-\bs{v}^{\prime})+g(\bs{v}^{\prime})}{\psi_2(\bs{v}^{\prime})>0,\; \bs{v}^{\prime}\in\mathbb{R}^{|F|}_+\setminus\{\bs{0}\}}, \\
&\mF_3=\Set*{U\le \nabla f(\bs{v}^{\prime})^{\top}(\bs{v}-\bs{v}^{\prime})+f(\bs{v}^{\prime})}{\psi_1(\bs{v}^{\prime})=0,\; \bs{v}^{\prime}\in\mathbb{R}^{|F|}_+\setminus\{\bs{0}\}}.
\end{aligned}
\end{equation}
\end{proposition}
\begin{proof}
Consider any $\bs{v}^{\prime}\in\mathbb{R}^{|F|}_+\setminus\{\bs{0}\}$ satisfying $\psi_1(\bs{v}^{\prime})>0$.
It is easy to see that the inequality $U\le \nabla f(\bs{v}^{\prime})\cdot(\bs{v}-\bs{v}^{\prime})+f(\bs{v}^{\prime})$ 
is a tangent inequality of $\chi_1$ with the tangent point at $(f(\bs{v}^{\prime}),\bs{v}^{\prime})$. 
Since $\psi_1(\bs{v}^{\prime})>0$, the distance between the hyper-plane 
$\Gamma=\Set*{(U,\bs{v})\in\mathbb{R}\times\mathbb{R}^{|F|}}
{U=\nabla f(\bs{v}^{\prime})\cdot(\bs{v}-\bs{v}^{\prime})+f(\bs{v}^{\prime})}$ and $\chi_2$
is positive. This indicates that $\Gamma$ is not tangent to $\chi_2$, 
which proves that the representation of $\mF_1$ in \eqref{eqn:F123} is valid.
Similarly, we can prove that the representation of $\mF_2$ in \eqref{eqn:F123} is valid.
For any $\bs{v}^{\prime}\in\mathbb{R}^{|F|}_+$ satisfying $\psi_1(\bs{v}^{\prime})=0$,
the hyper-plane $\Gamma$ is tangent to $\chi_2$ at the point $(g(\bs{v}^*),\bs{v}^*)$,
where $\bs{v}^*$ is the optimal solution of 
$\underset{\bs{v}}{\textrm{min}}\;\nabla f(\bs{v}^{\prime})\cdot(\bs{v}-\bs{v}^{\prime})+f(\bs{v}^{\prime})-g(\bs{v})$.
Therefore, $\Gamma$ is a common tangent plane of $\chi_1$ and $\chi_2$.
\end{proof}

The representation of $\conv(\chi)$ based on the tangent inequalities in $\mF_1$, $\mF_2$ and $\mF_3$ is
given in Theorem~\ref{thm:conv-hull-tang-plane}.
\begin{theorem}\label{thm:conv-hull-tang-plane}
Define the sets $\mathcal{V}_1$ and $\mathcal{V}_2$ as 
\begin{equation}\label{eqn:V1-V2}
\begin{aligned}
&\mathcal{V}_1=\Set*{\bs{v}^0\in\mathbb{R}^{|F|}_+\setminus\{\bs{0}\}}{\psi_1(\bs{v}_0)\ge 0}, \\
&\mathcal{V}_2=\Set*{\bs{v}^0\in\mathbb{R}^{|F|}_+\setminus\{\bs{0}\}}{\psi_2(\bs{v}_0)\ge 0}. 
\end{aligned}
\end{equation}
The convex hull $\emph{conv}(\chi)$ has the following representation:
\begin{equation}\label{eqn:convS=inter(v,U)}
\begin{aligned}
\emph{conv}(\chi)=\;&\mathcal{U}\cap\big(\cap_{\bs{v}^0\in\mathcal{V}_1}\Set*{(U,\bs{v})}{
U\le \nabla f(\bs{v}^0)^{\top}(\bs{v}-\bs{v}^0)+f(\bs{v}^0)}\big) \\
&\cap \left(\cap_{\bs{v}^0\in\mathcal{V}_2}\Set*{(U,\bs{v})}{
U\le \nabla g(\bs{v}^0)^{\top}(\bs{v}-\bs{v}^0)+g(\bs{v}^0)}\big)\right),
\end{aligned}
\end{equation}
where $\mathcal{U}=\Set*{(U,\bs{v})}{U\ge0,\;\bs{v}\ge\bs{0}}$.
\end{theorem}
\begin{proof}
Denote the set on the right side of \eqref{eqn:convS=inter(v,U)} as $\W$.
Clearly, we have $(0,\bs{0})\in\conv(\chi)$ and $(0,\bs{0})\in\W$.
It suffices to show that $\conv(\chi)\setminus\{(0,\bs{0})\}=\W\setminus\{(0,\bs{0})\}$.

We first show that $\textrm{conv}(\chi_1\cup\chi_2)\setminus\{(0,\bs{0})\}\subseteq \W\setminus\{(0,\bs{0})\}$. 
Note that for any point $\bs{v}^0\in\mathcal{V}_1$,
the tangent plane $U= \nabla f(\bs{v}^0)^{\top}(\bs{v}-\bs{v}^0)+f(\bs{v}^0)$
of $\chi_1$ at the point $(f(\bs{v}^0),\bs{v}^0)$ is also a tangent plane of $\textrm{conv}(\chi_1\cup\chi_2)$.
It indicates that the inequality $U\le \nabla f(\bs{v}^0)^{\top}(\bs{v}-\bs{v}^0)+f(\bs{v}^0)$
is a tangent inequality of $\textrm{conv}(\chi_1\cup\chi_2)$. 
Similarly, the inequality $U\le \nabla g(\bs{v}^0)^{\top}(\bs{v}-\bs{v}^0)+g(\bs{v}^0)$
is a tangent inequality of $\textrm{conv}(\chi_1\cup\chi_2)$ for any $\bs{v}^0\in\mathcal{V}_2$. 
Therefore, we have $\textrm{conv}(\chi_1\cup\chi_2)\setminus\{(0,\bs{0})\}\subseteq \W\setminus\{(0,\bs{0})\}$.

We now show that $\W\setminus\{(0,\bs{0})\}\subseteq\textrm{conv}(\chi_1\cup\chi_2)\setminus\{(0,\bs{0})\}$. 
Consider any point $(U^{\prime},\bs{v}^{\prime})\in \W\setminus\{(0,\bs{0})\}$.
We need to show that $(U^{\prime},\bs{v}^{\prime})\in\textrm{conv}(\chi_1\cup\chi_2)$. 
We prove it by contradiction. Suppose $(U^{\prime},\bs{v}^{\prime})$ is not in $\textrm{conv}(\chi_1\cup\chi_2)$.
Since $\textrm{conv}(\chi_1\cup\chi_2)$ is a closed convex set, by the separation principle 
(Theorem~11.1 in \citep{rockafellar1996_conv-analy}), 
there exists a plane $\Lambda=\Set*{(U,\bs{v})}{U-\bs{a}^{\top}\bs{v}=b}$, 
such that $U^{\prime}-\bs{a}^{\top}\bs{v}^{\prime}>b$ and $U-\bs{a}^{\top}\bs{v}\le b$ for all 
$(U,\bs{v})\in\textrm{conv}(\chi_1\cup\chi_2)$. Since $U-\bs{a}^{\top}\bs{v}\le b$ for all
$(U,\bs{v})\in\textrm{conv}(\chi_1\cup\chi_2)$, we can choose the parameter $b$ such that 
$b=\inf\{b^{\prime}:\;b^{\prime}\textrm{ satisfying }U-\bs{a}^{\top}\bs{v}\le b^\prime\;\;\forall (U,\bs{v})\in\textrm{conv}(\chi_1\cup\chi_2)\}$.
Since the set $\textrm{conv}(\chi_1\cup\chi_2)$ is closed, the infimum is attainable in the above set definition.
Therefore, the parameter $b$ can be chosen such that there exists a point $(U^0,\bs{v}^0)\in\textrm{conv}(\chi_1\cup\chi_2)$,
satisfying $U^0-\bs{a}^{\top}\bs{v}^0=b$, i.e., the point $(U^0,\bs{v}^0)$ is on the plane $\Lambda$.
We claim that there exists a point $(U^1,\bs{v}^1)\in (\mS_1\cup\mS_2)\cap\Lambda$, where
$\mS_1=\Set*{(U,\bs{v})\in\mathbb{R}\times\mathbb{R}^{|F|}_+\setminus\{(0,\bs{0})\}}
{U-f(\bs{v})=0}$ and
$\mS_2=\Set*{(U,\bs{v})\in\mathbb{R}\times\mathbb{R}^{|F|}_+\setminus\{(0,\bs{0})\}}
{U-g(\bs{v})=0}$ associated with the functions $f$ and $g$.
We will prove this claim at the end.
Without loss of generality, assume that $(U^1,\bs{v}^1)\in\mS_1\cap\Lambda$.
The proof for the case that $(U^1,\bs{v}^1)\in\mS_2\cap\Lambda$ is similar.
Since  $\mS_1$ is differentiable at $(U^1,\bs{v}^1)$,
the plane $\Lambda$ must be the tangent plane of $\mS_1$ at $(U^1,\bs{v}^1)$, which can be
written as: $U\le \nabla f(\bs{v}^1)^{\top}(\bs{v}-\bs{v}^1)+f(\bs{v}^1)$.
Since the plane $\Lambda$ separates the point $(U^{\prime},\bs{v}^{\prime})$ from the set $\textrm{conv}(\chi_1\cup\chi_2)$,
we must have $U^{\prime}> \nabla f(\bs{v}^1)^{\top}(\bs{v}^{\prime}-\bs{v}^1)+f(\bs{v}^1)$.
Since the inequality $U\le \nabla f(\bs{v}^1)^{\top}(\bs{v}-\bs{v}^1)+f(\bs{v}^1)$ is identical to 
$U-\bs{a}^\top\bs{v}\le b$, it follows that $U\le \nabla f(\bs{v}^1)^{\top}(\bs{v}-\bs{v}^1)+f(\bs{v}^1)$
for all $(U,\bs{v})\in\textrm{conv}(\chi_1\cup\chi_2)$. 
It indicates that $g(\bs{v})\le \max\{f(\bs{v}),\;g(\bs{v})\} \le \nabla f(\bs{v}^1)^{\top}(\bs{v}-\bs{v}^1)+f(\bs{v}^1)$ 
for all $\bs{v}$, and hence $\psi_1(\bs{v}^1)\ge 0$. It follows that $\bs{v}^1\in\mathcal{V}_1$ by definition.
Since $U^{\prime}> \nabla f(\bs{v}^1)^{\top}(\bs{v}^{\prime}-\bs{v}^1)+f(\bs{v}^1)$, 
it implies that at least one of the inequalities on the right side of \eqref{eqn:convS=inter(v,U)} 
has been violated by $(U^1,\bs{v}^1)$, which contradicts with $(U^{\prime},\bs{v}^{\prime})\in \textcolor{blue}{\W}$. 
Therefore, we have $\W\subseteq\textrm{conv}(\chi_1\cup\chi_2)$.

We now prove the claim that there exists a point $(U^1,\bs{v}^1)\in (\mS_1\cup\mS_2)\cap\Lambda$. 
Specifically, we need to show that if there exists a point $(U^0,\bs{v}^0)\in\textrm{conv}(\chi_1\cup\chi_2)\cap\Lambda$, 
there exists a point $(U^1,\bs{v}^1)\in(\mS_1\cup\mS_2)\cap\Lambda$. 
The point  $(U^0,\bs{v}^0)$ can be written as:
\begin{equation}\label{eqn:v0U0-conv-hull}
(U^0,\bs{v}^0)=\lambda(U^1,\bs{v}^1)+(1-\lambda)(U^2,\bs{v}^2),
\end{equation}
where $\lambda\in(0,1)$, $(U^1,\bs{v}^1)\in\chi_1$, and $(U^2,\bs{v}^2)\in\chi_2$.
Using \eqref{eqn:v0U0-conv-hull} and the fact that $(U^0,\bs{v}^0)\in\Lambda$, we have
\begin{equation}\label{eqn:sum_av-U=b}
\lambda(U^1-\bs{a}^{\top}\bs{v}^1)+(1-\lambda)(U^2-\bs{a}^{\top}\bs{v}^2)=b.
\end{equation}
Since we also have $U^1-\bs{a}^{\top}\bs{v}^1\le b$ and $U^2-\bs{a}^{\top}\bs{v}^2\le b$,
it combined with \eqref{eqn:sum_av-U=b} implies that $U^1-\bs{a}^{\top}\bs{v}^1=b$
and $U^2-\bs{a}^{\top}\bs{v}^2=b$. 
Therefore, we have $(\bs{v}^1,U^1)\in \mS_1\cap\Lambda\subseteq(\mS_1\cup\mS_2)\cap\Lambda$,
which concludes the proof of the claim.
\end{proof}

The representation of $\conv(\chi)$ given in Theorem~\ref{thm:conv-hull-tang-plane}
provides a computational framework for generating the tangent inequalities of $\conv(\chi)$ in an algorithm for solving \eqref{opt:RFL}.
The framework takes the current optimal solution $\bs{v}_0$ 
(only the $\bs{v}$-component of the solution matters) as the input, and generates tangent inequalities of $\conv(\chi)$ based on this point.
For a given $\bs{v}_0$, we can generate the following gradient based inequalities (\textbf{G-cuts}):
\begin{align}
&U\le \nabla f(\bs{v}_0)^{\top}(\bs{v}-\bs{v}_0)+f(\bs{v}_0), \label{eqn:ineq-grad1} \\
&U\le \nabla g(\bs{v}_0)^{\top}(\bs{v}-\bs{v}_0)+g(\bs{v}_0). \label{eqn:ineq-grad2}
\end{align}
The inequality \eqref{eqn:ineq-grad1} is a tangent inequality of $\conv(\chi)$ if and only if $\psi_1(\bs{v}_0)\ge 0$,
and the inequality \eqref{eqn:ineq-grad2} is a tangent inequality of $\conv(\chi)$ if and only if $\psi_2(\bs{v}_0)\ge 0$. The third type of inequality can be generated using a disjunctive formulation. This formulation is given in the next subsection.

\subsection{Convexification using a Disjunctive Formulation}
\label{sec:conv-2}
The set $\conv(\chi)$ can alternatively be represented based on the lift-and-project technique 
that is widely used in the research of mixed integer  programming~\citep{balas1998,stub2002_BC-0-1-mixed-conv-prog}. 
A tangent inequality is induced by 
a point $(U_0,\bs{v}_0)$ outside $\conv(\chi)$. We can construct a convex optimization
problem to generate a tangent plane of $\conv(\chi)$ that separates $(U_0,\bs{v}_0)$ from $\conv(\chi)$.
This is given in the following Proposition~\ref{prop:disj-prog}.
\begin{proposition}\label{prop:disj-prog}
Let $(U_0,\bs{v}_0)$ be any point outside $\emph{conv}(\chi)$.
Solving the following convex optimization problem to get an optimal solution
$(U^*,\bs{v}^*,U^*_1,\bs{v}^*_1,U^*_2,\bs{v}^*_2)$:
\begin{equation}\label{opt:min-dist}
\begin{aligned}
&\emph{min}\;\; \|\bs{v}_0-\bs{v}\|^2+(U_0-U)^2  \\
&\emph{ s.t. }\;\; \bs{v} = \bs{v}_1+\bs{v}_2,\quad U=U_1+U_2, \\
&\qquad \;\; U_1\le f(\bs{v}_1),\quad U_2\le g(\bs{v}_2), \\
&\qquad \;\; U,U_1,U_2\ge 0,\; \bs{v},\bs{v}_1,\bs{v}_2\ge 0.
\end{aligned}
\end{equation}
If $\norm{\bs{v}^*_1}>0$ and $\norm{\bs{v}^*_2}>0$,
the following inequality is a tangent inequality of $\emph{conv}(\chi)$:
\begin{equation}\label{eqn:lift-proj-facet}
U\le -\frac{1}{U_0-U^*}(\bs{v}_0-\bs{v}^*)^{\top} (\bs{v}-\bs{v}^*)+U^*,
\end{equation}
\end{proposition}
\begin{proof}
We first show that the convex hull $\conv(\chi)=\textrm{conv}(\chi_1\cup\chi_2)$ can be represented in the following form:
\begin{equation}\label{eqn:conv-hull-rep}
\textrm{conv}(\chi_1\cup\chi_2)=\Set*{(U,\bs{v})}{\arraycolsep=5pt\def\arraystretch{1.5}
											\begin{array}{ll}
											U=U_1+U_2, & \bs{v}=\bs{v}_1+\bs{v}_2, \\
											0\le U_1\le f(\bs{v}_1), & \bs{v}_1\ge \bs{0}, \\
											0\le U_2\le g(\bs{v}_2), & \bs{v}_2\ge \bs{0}
											 \end{array}
}.
\end{equation}
Let the set on the right side of \eqref{eqn:conv-hull-rep} be $\chi$. To show that 
$\textrm{conv}(\chi_1\cup\chi_2)\subseteq \chi$, we let $(U,\bs{v})$ be any point in $\textrm{conv}(\chi_1\cup\chi_2)$
and show that $(U,\bs{v})\in \chi$. 
Since $(U,\bs{v})\in\textrm{conv}(\chi_1\cup\chi_2)$, 
there exist $\lambda\in [0,1]$, $(U^{\prime}_1,\bs{v}^{\prime}_1)\in\chi_1$ and $(U^{\prime}_2,\bs{v}^{\prime}_2)\in\chi_2$ such that 
\begin{displaymath}
\bs{v}=\lambda\bs{v}^{\prime}_1+(1-\lambda)\bs{v}^{\prime}_2,\quad U=\lambda U^{\prime}_1+(1-\lambda)U^{\prime}_2.
\end{displaymath}
Let $U_1=\lambda U^{\prime}_1$, $\bs{v}_1=\lambda\bs{v}_1$, $U_2=\lambda U^{\prime}_2$, and $\bs{v}_2=\lambda\bs{v}_2$.
Since $0\le U^{\prime}_1\le f(\bs{v}^\prime_1)$, we have $0\le U_1=\lambda U^{\prime}_1\le \lambda f(\bs{v}^{\prime}_1)=f(\lambda \bs{v}^{\prime}_1)=f(\bs{v}_1)$
which is based on the observation that the function $f$ is scale-invariant. Similarly, we can show that $0\le U_2\le g(\bs{v}_2)$.
Therefore, we have shown that $(U,\bs{v})\in \chi$. To show that $\chi\subseteq \textrm{conv}(\chi_1\cup\chi_2)$,
we let $(U,\bs{v})$ be any point in $\chi$, and show $(U,\bs{v})\in\textrm{conv}(\chi_1\cup\chi_2)$.
Since $(U,\bs{v})\in\chi$, there exist $U^{\prime}_1,U^{\prime}_2,\bs{v}^{\prime}_1,\bs{v}^{\prime}_2$ such that
\begin{displaymath}
\begin{aligned}
&U=U^{\prime}_1+U^{\prime}_2, \qquad \bs{v}=\bs{v}^{\prime}_1+\bs{v}^{\prime}_2, \\
&0\le U^{\prime}_1\le f(\bs{v}_1), \qquad \bs{v}^{\prime}_1\ge \bs{0}, \\
&0\le U^{\prime}_2\le g(\bs{v}_2), \qquad \bs{v}^{\prime}_2\ge \bs{0}.
\end{aligned}
\end{displaymath}
Let $U_1=2U^{\prime}_1$, $\bs{v}_1=2\bs{v}^{\prime}_1$, $U_2=2U^{\prime}_2$, and $\bs{v}_2=2\bs{v}^{\prime}_2$,
we have $U=(U_1+U_2)/2$, $\bs{v}=(\bs{v}_1+\bs{v}_2)/2$, $0\le U_1=2U^{\prime}_1\le2f(\bs{v}^{\prime}_1)=f(2\bs{v}^{\prime}_1)=f(\bs{v}_1)$,
and similarly $0\le U_2\le g(\bs{v}_2)$. Therefore, we have shown that $(U,\bs{v})\in\textrm{conv}(\chi_1\cup\chi_2)$.
The representation \eqref{eqn:conv-hull-rep} is proved.

Notice that the optimal solution $(U^*,\bs{v}^*)$ of the convex optimization problem \eqref{opt:min-dist}
is the point in $\textrm{conv}(\chi_1\cup\chi_2)$ that has the shortest distance (measured by the $\ell_2$-norm) with respect to the point 
$(U^0,\bs{v}^0)$.
Consider the tangent plane $\Gamma$ of $\textrm{conv}(\chi_1\cup\chi_2)$ that passes the point $(U^*,\bs{v}^*)$.
Since $\norm{(U^*_1,\bs{v}^*_1)}>0$ and $\norm{(U^*_2,\bs{v}^*_2)}>0$,  $\Gamma$ must be the tangent plane
of $\chi_1$ at $(U^*_1,\bs{v}^*_1)$ and the tangent plane of $\chi_2$ at $(U^*_2,\bs{v}^*_2)$.
By the basic result from analytical geometry, we know that the normal vector of $\Gamma$ is given by
$(U^*-U_0,\bs{v}^*-\bs{v}_0)$ up to a scale factor, and hence the plane $\Gamma$ can be written as: 
\begin{displaymath}
(U^0-U^*)(U-U^*)+(\bs{v}^0-\bs{v}^*)^{\top}\cdot(\bs{v}-\bs{v}^*)=0,
\end{displaymath}
which implies that \eqref{eqn:lift-proj-facet} is a tangent inequality of $\textrm{conv}(\chi_1\cup\chi_2)$.
\end{proof}
The tangent inequality generated using the lift-and-project technique in Proposition~\ref{prop:disj-prog}
depends on a point not in $\conv(\chi)$.  Theorem~\ref{thm:conv-lift-and-project} shows that 
all tangent inequalities in $\mF$ can be generated using Proposition~\ref{prop:disj-prog} to construct $\conv(\chi)$.
\begin{theorem}\label{thm:conv-lift-and-project}
Let $\mtc{U}=\Set*{(U,\bs{v})}{U\ge0,\;\bs{v}\ge\bs{0}}$ and $\mtc{S}=\mtc{U}\setminus\emph{conv}(\chi)$.
The convex hull $\emph{conv}(\chi)$ can be represented as:
\begin{equation}\label{eqn:rep-conv-disj}
\emph{conv}(\chi) = \mtc{U}\cap\left(\cap_{(U_0,\bs{v}_0)\in\mtc{S}}
						 \Set*{(U,\bs{v})}{U\le \alpha(U_0,\bs{v}_0)^\top\bs{v}+\beta(U_0,\bs{v}_0) }
						  \right),
\end{equation}
where $\alpha(U_0,\bs{v}_0)$ and $\beta(U_0,\bs{v}_0)$ are coefficients that are functions of $(U_0,\bs{v}_0)$
determined in the following steps:
\begin{enumerate}
	\item Solve the projection problem \eqref{opt:min-dist} to get an optimal solution 
		$(U^*,U^*_1,\bs{v}^*_1,U^*_2,\bs{v}^*_2)$;
	\item If $\norm{(U^*_1,\bs{v}^*_1)}>0$ and $\norm{(U^*_2,\bs{v}^*_2)}>0$, let
		$\alpha(U_0,\bs{v}_0)=-\frac{1}{U_0-U^*}(\bs{v}_0-\bs{v}^*)$ and 
		$\beta(U_0,\bs{v}_0)=\frac{1}{U_0-U^*}(\bs{v}_0-\bs{v}^*)^\top\bs{v}^*+U^*$;
	\item If $\norm{(U^*_1,\bs{v}^*_1)}>0$ but $\norm{(U^*_2,\bs{v}^*_2)}=0$,
		let $\alpha(U_0,\bs{v}_0)=\nabla f(\bs{v}_0)$ 
		and $\beta(U_0,\bs{v}_0)=-\nabla f(\bs{v}_0)^\top\bs{v}_0+f(\bs{v}_0)$;
	\item If $\norm{(U^*_1,\bs{v}^*_1)}=0$ but $\norm{(U^*_2,\bs{v}^*_2)}>0$,
		let $\alpha(U_0,\bs{v}_0)=\nabla g(\bs{v}_0)$ 
		and $\beta(U_0,\bs{v}_0)=-\nabla g(\bs{v}_0)^\top\bs{v}_0+g(\bs{v}_0)$.
\end{enumerate}
\end{theorem}
\begin{proof}
Let $V$ be the set on the right side of \eqref{eqn:rep-conv-disj}.
It is proved in Proposition~\ref{prop:disj-prog} that $\conv(\chi)\subseteq V$.
It suffices to show that $V\subseteq\conv(\chi)$. We prove it by contradiction.
Suppose there exists a point $(U_0,\bs{v}_0)\in V\setminus\conv(\chi)$.
Clearly, we have $(U_0,\bs{v}_0)\in\mtc{S}$. By solving \eqref{opt:min-dist} at $(U_0,\bs{v}_0)$,
we construct a tangent plane 
\begin{equation}
U=\alpha(U_0,\bs{v}_0)^\top\bs{v}+\beta(U_0,\bs{v}_0),
\end{equation}
where the coefficients $\alpha(U_0,\bs{v}_0)$ and $\beta(U_0,\bs{v}_0)$
are determined as above based on three different cases. In the first case,
$\bs{v}^*=\bs{v}^*_1+\bs{v}^*_2$ for non-zero $\bs{v}^*_1$ and $\bs{v}^*_2$,
and hence the tangent plane has the form \eqref{eqn:lift-proj-facet}.
In the second and third cases, either $\bs{v}^*_1$ or $\bs{v}^*_2$ is zero,
and the tangent plane is in the form of \eqref{eqn:ineq-grad1} and \eqref{eqn:ineq-grad1},
respectively. In any case, the tangent plane $U=\alpha(U_0,\bs{v}_0)^\top\bs{v}+\beta(U_0,\bs{v}_0)$
of $\conv(\chi)$ strictly separates $(U_0,\bs{v}_0)$ from $\conv(\chi)$, which implies that the 
inequality $U\le\alpha(U_0,\bs{v}_0)^\top\bs{v}+\beta(U_0,\bs{v}_0)$ is violated by the point $(U,\bs{v})=(U_0,\bs{v}_0)$.
This contradicts with $(U_0,\bs{v}_0)\in V$.
\end{proof}

The lift-and-project technique in Proposition~\ref{prop:disj-prog} can be used as a common approach to generate tangent inequalities 
in $\mF_1$, $\mF_2$ and $\mF_3$ if for a given $\bs{v}_0$ the point $(\hat{\bs{\beta}}^{\top}\bs{v}_0,\bs{v}_0)$ is outside  $\conv(\chi)$.
Let $(U^*,\bs{v}^*)$ be the optimal solution of the following convex program:
\begin{equation}\label{opt:min-dist-2}
\begin{aligned}
&\textrm{min}\;\; \|\bs{v}_0-\bs{v}\|^2+(\hat{\bs{\beta}}^{\top}\bs{v}_0-U)^2  \\
&\textrm{ s.t. }\;\; \bs{v} = \bs{v}_1+\bs{v}_2,\quad U=U_1+U_2, \\
&\qquad \;\; U_1\le f(\bs{v}_1),\quad U_2\le g(\bs{v}_2), \\
&\qquad \;\; U,U_1,U_2\ge 0,\; \bs{v},\bs{v}_1,\bs{v}_2\ge 0.
\end{aligned}
\end{equation}
If $\norm{\bs{v}^*_1}>0$ and $\norm{\bs{v}^*_2}>0$,
we can add the following lift-and-project inequality (\textbf{P-cut}) to convexify \eqref{eqn:nonconv-constr}:
\begin{equation}\label{eqn:P-cut}
U\le -\frac{1}{\hat{\bs{\beta}}^{\top}\bs{v}_0-U^*}(\bs{v}_0-\bs{v}^*)^{\top}(\bs{v}-\bs{v}^*)+U^*.
\end{equation}
Note that the convex optimization problem \eqref{opt:min-dist-2} can be reformulated as
a convex quadratic-constraint-quadratic-programming problem, which can be solved using Gurobi~\citep{gurobi}.

Based on the tangent cuts developed in this section, we can establish a cutting-plane algorithm 
for solving \eqref{opt:RFL-relax} iteratively. The pseudo code is given in Algorithm~\ref{alg:cut-plane}.
\begin{algorithm}
	{\footnotesize
	\caption{\footnotesize A cutting-plane algorithm for solving \eqref{opt:RFL-relax}.}
	\label{alg:cut-plane}
	\begin{algorithmic}[1]
	\State{Input: a threshold value $\epsilon$ of error tolerance.}
	\State{Solving the relaxation problem of \eqref{opt:RFL-relax} by replacing the constraints 
	 $U^{ij}\le\psi^{ij}(\bs{v}^{ij})$ with $U^{ij}\le \hat{\beta}^\top\bs{v}^{ij}$. 
	 Let $(U^*,\bs{y}^*,\bs{v}^*)$ be an optimal solution of the initial relaxation problem.}
	 \State{Let $(U^\prime,\bs{y}^\prime,\bs{v}^\prime)\gets(U^*,\bs{y}^*,\bs{v}^*)$ 
	 	and $\delta=\infty$. }
	\While{$\delta>\epsilon$}
		\State{Derive tangent inequalities for each $U^{ij}\le\psi^{ij}(\bs{v}^{ij})$ induced by the current
			optimal solution $\bs{v}^*$ in the following order:}
			\begin{itemize}
				\item[-] If $\psi_1(\bs{v}^{*ij})\ge 0$, generates a $\mF_1$-type G-cut
				$U^{ij}\le \nabla f(\bs{v}^{*ij})^\top(\bs{v}^{ij}-\bs{v}^{*ij}) + f(\bs{v}^{*ij})$ and go to Line~\ref{lin:add-cut};
				\item[-] If $\psi_2(\bs{v}^{*ij})\ge 0$, generates a $\mF_2$-type G-cut cut
				$U^{ij}\le \nabla g(\bs{v}^{*ij})^\top(\bs{v}^{ij}-\bs{v}^{*ij}) + g(\bs{v}^{*ij})$ and go to Line~\ref{lin:add-cut};
				\item[-] Solve the convex program (32) with $\bs{v}_0\gets\bs{v}^{*ij}$ to generate a P-cut (33).
			\end{itemize}
		\State{Update the relaxation problem of \eqref{opt:RFL-relax} by adding the newly generated tangent 				inequalities. } \label{lin:add-cut}
		\State{Solve the update relaxation problem to get an optimal solution $(U^*,\bs{y}^*,\bs{v}^*)$.}
		\State{Set $\delta\gets\norm{\bs{v}^*-\bs{v}^\prime}_{\infty}$ and 
			$(U^\prime,\bs{y}^\prime,\bs{v}^\prime)\gets(U^*,\bs{y}^*,\bs{v}^*)$.}
	\EndWhile
	\State{Return $\bs{y}^*$ as an optimal decision vector of \eqref{opt:RSP-relax}.}
	\end{algorithmic}
	}
\end{algorithm}

\begin{remark}
In practice, the implementation of Algorithm~\ref{alg:cut-plane} will be run 
for a limited amount of computational time. The best solution will be returned
in the case that the algorithm does not terminate (or converge) in the given 
amount of time. See Table~\ref{tab:determ-dmd} for the numerical results 
lead by the practical implementation of the algorithm. We also note that 
although Algorithm~\ref{alg:cut-plane} for solving \eqref{opt:RFL-relax} 
can terminate faster and often lead to a better solution than directly solving
\eqref{opt:RFL-determ} in a branch-and-cut framework, 
\eqref{opt:RFL-determ} is an exact reformulation of \eqref{opt:RFL},
whereas \eqref{opt:RFL-relax} is an approximation problem. 
\end{remark}

\section{Computational Experiments}
\label{sec:comp-exp}
In this section we discuss computational performance of the proposed approach and its implications on location decisions using randomly generated problems. Problem instances generated from a Covid-19 case study data are discussed in the next section.
\subsection{Numerical Instance Generation}
\label{sec:gen-inst}
We generated 41 \eqref{opt:RFL} instances to test the computational performance of solving the 
\eqref{opt:RFL-determ} using the developed techniques. The instances are labeled as FL0, FL1, \ldots, FL41.
The FL0$\sim$FL25 instances are small and mid-size, and FL26$\sim$FL40 are large instances
in terms of customer sites, candidate locations and total budget.  
An instance is determined by the following parameters: 
number of customer locations $|S|$, number of candidate service center locations $|F|$, 
the total budget $B$ for establishing the facilities,
the capacity $C_j$ of each service center, the demand $D_i$ for each customer site,
and all the parameters for determining
the ambiguity set \eqref{eqn:Puy-reform} for all $i\in S$ and $j\in F$.

We now describe the numerical instance generation. 
The number of customer sites $|S|$ is given in the second column of Table~\ref{tab:determ-dmd}. 
The customer sites are points located in a $15\times 15$ two-dimensional square. 
The two coordinates of each customer site are generated using a uniform random variable in the range $[0,15]^2$. 
Every customer site is also a candidate service center location, i.e., $F=S$. 
The parameters $\bs{c}$ that represent the extra gain in establishing service centers
in the \eqref{opt:RFL} model are set to zero in all the numerical instances. Therefore, the instances only consider the total expected utility gained by customers.  
The cost of establishing each service center is 1, i.e., $b_j=1$ for all $j\in F$ in \eqref{opt:RFL}. 
The total budget is given in the third column of Table~\ref{tab:determ-dmd}.
For every $j\in F$, the capacity $C_j$ is generated randomly from the interval $[100,180]$.
To define the parameters in \eqref{eqn:Puy-reform}, we first define an effective distance $L_0=5$,
and define an effective set $F_i$ of service centers for each $i\in S$ such that 
$F_i=\Set*{j\in F}{\|\bs{x}^j-\bs{x}^i\|_2\le L_0}$, where $\bs{x}^i$ is the coordinate vector of the
customer site $i\in S$.  The parameters $\hat{\bs{\beta}}^{ij}$ are set as follows:
\begin{equation}
\hat{\beta}^{ij}_k = \left\{\def\arraystretch{1.2}
			 \begin{array}{ll}  
			  10\times \left(1-\|\bs{x}^i-\bs{x}^j\|_2/L_0\right)  & \qquad \textrm{if } j\in F_i\textrm{ and }k=j  \\
			  1-\|\bs{x}^i-\bs{x}^k\|_2/L_0 & \qquad \textrm{if } j\in F_i\textrm{ and }k\neq j  \\
			  0 & \qquad \textrm{if } j\in F\setminus F_i,\;\forall k\in F.
			  \end{array}
			  \right.
\end{equation}
Thus, the parameters reflect inverse proportionality to utility with respect to distance. The covariance matrix $\widehat{\bs{\Sigma}}^{ij}$ (for all $i\in S,\;j\in F$) is set to be $\widehat{\bs{\Sigma}}^{ij}=\bs{Q}^{ij\top}\bs{Q}^{ij}$,
where $\bs{Q}^{ij}$ is a $|F|\times |F|$ matrix with each entry randomly generated from $[0,1]$.
The matrix $\bs{A}^{ij}$ (for all $i\in S,\;j\in F$) is set to be $\bs{A}^{ij}=\bs{I}_{|F|}+0.3\;\bs{Q}^{ij\top}\bs{Q}^{ij}$,
where $\bs{I}_{|F|}$ is the $|F| \times |F|$ identity matrix. We set $\gamma^{ij}_1=0.05$, 
$\gamma^{ij}_2=0.2$ and $b^{ij}=0.2$ for all $i\in S$, $j\in F$.

\subsection{Computational Performance of Solving \eqref{opt:RFL} Instances}
\label{sec:comp-determ-dmd}
We conducted experiments to test the computational performance of solving 
\eqref{opt:RFL} instances with two different approaches. The first approach
is to reformulate it into \eqref{opt:RFL-II-2} which is a MISOCP in a lifted space,
and solve it directly using a MISOCP solver (e.g., Gurobi). The second approach is try to solve
the relaxed problem \eqref{opt:RSP-relax} using the cutting-plane algorithm (Algorithm~\ref{alg:cut-plane}).
At each iteration of Algorithm~\ref{alg:cut-plane}, the relaxation problem is a MILP which is solved using Gurobi.
The threshold value $\epsilon$ for convergence in Algorithm~\ref{alg:cut-plane} is set to be $0.001$ for all experiments.
Table~\ref{tab:determ-dmd} compares the computational performance of solving 30 
\eqref{opt:RFL-determ} instances by using the lifted MISOCP reformulation 
 versus using the cutting-plane algorithm (Algorithm~\ref{alg:cut-plane}).
 The comparison shows that the cutting-plane algorithm identifies a better solution
 in 25 instances, the lifted MISOCP identifies a better solution in only 1 instance,
 and the two approaches are in tie in 4 instances. In particular, the lifted MISOCP becomes
 intractable for the large instances such as FL26$\sim$FL30, for which the objective returned
 by solving the lifted MISOCP is in magnitude smaller than that from the cutting-plane algorithm.   
 This indicates that the cutting-plane method is more effective in practice. We also found that
 the relaxed MILP in the cutting-plane algorithm is much easier to solve than the lifted MISOCP
 due to that the former has a much smaller model size and no SOC constraints. 
 This is also reflected in the number of nodes and solver cuts. Notice that the cutting-plane
 algorithm converges in 12 instances within the 4-hour time limit, indicating that it may converge
 in more instances if more computational time is given. A considerable amount of tangent cuts 
 have been generated to strengthen the MILP relaxation in the algorithm.

To get a better understanding of which types of tangent cuts have an essential contribution
to strengthen the MILP relaxation in the cutting-plane algorithm, we investigate two approaches
of adding the cuts. In the first approach, we only add the $\mtc{F}_1$ and $\mtc{F}_2$ cuts to 
the MILP, whereas in the second approach, we add $\mtc{F}_1$, $\mtc{F}_2$ and $\mtc{F}_3$ cuts. 
The computational comparison of the two approaches is summarized in Table~\ref{tab:cut12_vs_cut123}.
The results show that there is not much difference in the best objective identified by the two approaches. 
Only for the instances FL15, FL26 and FL29, the second approach identifies a better objective value,
while for the instances FL25 and LF30, the first approach identifies a better objective value. 
We have also found that in all the calculation, not a single $\mtc{F}_2$ cut has been generated. 
This indicates that for all the numerical instances generated, the first SOC function $f^{ij}$ plays
a dominate role in determining $\max\{f^{ij}(\bs{v}),\;g^{ij}(\bs{v})\}$. Based on our experience
of generating numerical instances, it is rather difficult to create an instance such that the two SOC
functions $f^{ij}$ and $g^{ij}$ are non-trivially competing to determine the value of 
$\max\{f^{ij}(\bs{v}),\;g^{ij}(\bs{v})\}$ without a very deliberate
and artificial tuning of the coefficients in $f^{ij}$ and $g^{ij}$. This indicates that in a naturally generated
instance, it is very likely that one of the SOC functions will dominate the other 
and hence the cutting-plane algorithm may identify an optimal solution of the original
problem. Although adding $\mtc{F}_3$ cuts lead to a marginal improvement on the 
objective value and the relaxation objective value, solving the cut-generation problem \eqref{opt:min-dist-2}
can take extra time. This is reflected in the CPU time comparison in Table~\ref{tab:cut12_vs_cut123}.
With 4-hour computational time limit, the approach of adding $\mtc{F}_1$ and $\mtc{F}_2$ cuts leads 
to 17 convergent instances, whereas there are only 10 convergent instances with the approach 
of adding $\mtc{F}_1$, $\mtc{F}_2$ and $\mtc{F}_3$ cuts.

\begin{table}
\centering
\begin{threeparttable}
{
\tiny
\caption{\footnotesize
Numerical results of solving 30 \eqref{opt:RFL} instances using the lifted MISOCP reformulation 
\eqref{opt:RFL-determ} versus using the relaxation of \eqref{opt:RFL-relax}  with $\mtc{F}_1$-, $\mtc{F}_2$- and $\mtc{F}_3$-type cuts developed in Section~\ref{sec:conv-1} and \ref{sec:conv-2}. 
For the two approaches, the computational time limit is set to be 4 hours (14400 seconds).
For the relaxed reformulation, 4 hours is the time limit of running Algorithm~\ref{alg:cut-plane}
and we also set the time limit of solving each relaxation MILP to be 30 mins. The symbol `-' in the column of `CPU(s)' indicates that the 4-hour computational time limit has been reached.  
For the two approaches, the column `nodes' is the number of nodes explored by the solver in the branch-and-bound process, the column `obj' is the best objective value found within the time limit,
and the column `solverCuts' is the number of cuts generated automatically by the solver. 
For the relaxed reformulation, the column `iters' is the number of iterations in Algorithm~\ref{alg:cut-plane}
for updating the relaxation problem,
the column `convg' indicates whether Algorithm~\ref{alg:cut-plane} converges or not within the 4 hours time limit,
the `cutRelaxObj' is the objective function of the relaxation problem in the last iteration of Algorithm~\ref{alg:cut-plane},
and the column `devCuts' (developed cuts) is the number of cuts generated 
using methods from Sections~\ref{sec:conv-1} and \ref{sec:conv-2}.
}\label{tab:determ-dmd}
\center
\setlength\tabcolsep{1.5pt}
\begin{tabular}{ccc|ccccc|ccccccccc}
\hline\hline
\multicolumn{3}{c}{Instance} & \multicolumn{5}{c}{lifted reformulation \eqref{opt:RFL-determ}} &\multicolumn{9}{c}{relaxed reformulation with $\mtc{F}_1$-, $\mtc{F}_2$- and $\mtc{F}_3$-type cuts} \\
\hline
ID	&	$|S|$	&	$B$	&	CPU(s)	&	obj	&	gap(\%)	&	nodes	&	solverCuts	&	CPU(s)	&	iters	&	convg	&	obj	&	cutRelaxObj	&	MipGap(\%)	&	nodes	&	solverCuts	&	devCuts	\\
\hline																																	
FL1	&	10	&	5	&	4	&	4248	&	0	&	92	&	62	&	5	&	5	&	y	&	4248	&	4312	&	0	&	1	&	21	&	40	\\
FL2	&	40	&	5	&	1029	&	6499	&	0	&	3482	&	394	&	175	&	5	&	y	&	6499	&	6499	&	0	&	388	&	87	&	99	\\
FL3	&	40	&	10	&	-	&	11862	&	15.96	&	4267	&	724	&	-	&	10	&	n	&	11926	&	12030	&	0	&	335602	&	0	&	393	\\
FL4	&	40	&	15	&	-	&	14217	&	15.19	&	14317	&	443	&	-	&	28	&	n	&	14250	&	14447	&	0	&	6271	&	310	&	1133	\\
FL5	&	60	&	5	&	4292	&	7094	&	0	&	2943	&	344	&	271	&	5	&	y	&	7094	&	7094	&	0	&	319	&	96	&	92	\\
FL6	&	60	&	10	&	-	&	13868	&	11.34	&	7271	&	615	&	-	&	20	&	n	&	13971	&	14033	&	0	&	7694	&	231	&	835	\\
FL7	&	60	&	15	&	-	&	18943	&	21.48	&	1788	&	728	&	-	&	9	&	n	&	19931	&	20163	&	0	&	16681	&	556	&	544	\\
FL8	&	80	&	5	&	4065	&	7409	&	0	&	3307	&	252	&	262	&	4	&	y	&	7409	&	7456	&	0	&	409	&	75	&	64	\\
FL9	&	80	&	10	&	-	&	14287	&	12.23	&	2882	&	585	&	5089	&	5	&	y	&	14404	&	14561	&	0	&	4200	&	206	&	188	\\
FL10	&	80	&	15	&	-	&	19928	&	19.77	&	1405	&	565	&	-	&	9	&	n	&	20281	&	20422	&	0	&	11280	&	586	&	551	\\
FL11	&	100	&	5	&	12815	&	6935	&	0	&	7336	&	540	&	1411	&	8	&	y	&	6905	&	6991	&	0	&	487	&	144	&	144	\\
FL12	&	100	&	10	&	-	&	13617	&	18.76	&	4311	&	545	&	-	&	17	&	n	&	14237	&	14349	&	0	&	3768	&	211	&	752	\\
FL13	&	100	&	15	&	-	&	19449	&	21.93	&	1572	&	686	&	10178	&	6	&	y	&	20410	&	20484	&	1.92	&	10380	&	519	&	347	\\
FL14	&	200	&	5	&	-	&	6957	&	19.93	&	2109	&	355	&	2587	&	5	&	y	&	7377	&	7498	&	0	&	558	&	209	&	93	\\
FL15	&	200	&	10	&	-	&	13384	&	21.67	&	2145	&	571	&	-	&	9	&	n	&	14122	&	14322	&	0	&	12014	&	298	&	356	\\
FL16	&	200	&	15	&	-	&	19882	&	21.75	&	2214	&	972	&	-	&	8	&	n	&	20850	&	21050	&	2.08	&	5985	&	503	&	478	\\
FL17	&	300	&	5	&	-	&	7459	&	17.84	&	3018	&	282	&	12848	&	9	&	y	&	7532	&	7571	&	1.35	&	2767	&	286	&	184	\\
FL18	&	300	&	10	&	-	&	14465	&	24.85	&	1373	&	975	&	9009	&	5	&	y	&	15397	&	15482	&	5.33	&	5869	&	304	&	193	\\
FL19	&	300	&	15	&	-	&	20517	&	31.53	&	1360	&	1057	&	-	&	8	&	n	&	22776	&	22949	&	3.96	&	2937	&	282	&	499	\\
FL20	&	400	&	5	&	-	&	7539	&	22.01	&	3228	&	496	&	-	&	8	&	y	&	7619	&	7695	&	7.66	&	5404	&	239	&	150	\\
FL21	&	400	&	10	&	-	&	14863	&	20.93	&	3364	&	788	&	-	&	8	&	n	&	15448	&	15597	&	6.54	&	2863	&	443	&	320	\\
FL22	&	400	&	15	&	-	&	21244	&	25.72	&	3248	&	1003	&	-	&	8	&	n	&	22743	&	22982	&	6.7	&	2928	&	552	&	497	\\
FL23	&	500	&	5	&	-	&	7645	&	17.06	&	2854	&	558	&	-	&	8	&	n	&	7674	&	7719	&	7.9	&	3484	&	256	&	163	\\
FL24	&	500	&	10	&	-	&	14617	&	24.94	&	3093	&	498	&	-	&	8	&	n	&	15268	&	15483	&	6.86	&	1731	&	385	&	309	\\
FL25	&	500	&	15	&	-	&	21619	&	25.66	&	3099	&	1012	&	-	&	8	&	n	&	22703	&	22932	&	6.59	&	2856	&	499	&	462	\\
FL26	&	1000	&	100	&	-	&	7715	&	na	&	1	&	2502	&	-	&	8	&	n	&	148060	&	148821	&	14.76	&	1	&	640	&	3116	\\
FL27	&	1000	&	200	&	-	&	14750	&	na	&	1	&	3173	&	-	&	8	&	n	&	244054	&	246025	&	31.52	&	0	&	0	&	4057	\\
FL28	&	1000	&	300	&	-	&	2876	&	na	&	1	&	1111	&	12604	&	7	&	y	&	317858	&	317910	&	60.2	&	1	&	830	&	5740	\\
FL29	&	1500	&	100	&	-	&	0	&	na	&	1	&	2824	&	-	&	8	&	n	&	143549	&	145698	&	14.98	&	1	&	284	&	3075	\\
FL30	&	1500	&	200	&	-	&	8750	&	na	&	1	&	4558	&	-	&	8	&	n	&	264751	&	268994	&	20.23	&	1	&	0	&	5148	\\
\hline
\end{tabular}
}
\end{threeparttable}
\end{table}

\begin{table}
\centering
\begin{threeparttable}
{
\tiny
\caption{\footnotesize 
Comparison of the computational performance between the case of
only adding $\mtc{F}_1$ and $\mtc{F}_2$ cuts versus adding $\mtc{F}_1$, $\mtc{F}_2$ and $\mtc{F}_3$ cuts
developed in Sections~\ref{sec:conv-1} and \ref{sec:conv-2}. Notice that the number of $\mtc{F}_2$ cuts
is zero for both cases (i.e., no $\mtc{F}_2$ cuts are generated), and hence there is no corresponding 
column for $\mtc{F}_2$ cuts. The symbol `-' in the column of `CPU(s)' indicates that the 4-hour computational time limit has been reached.
}\label{tab:cut12_vs_cut123}
\center
\setlength\tabcolsep{2.0pt}
\begin{tabular}{ccc|cccccc|ccccccc}
\hline\hline
\multicolumn{3}{c}{Instance} & \multicolumn{6}{c}{Adding $\mtc{F}_1$ and $\mtc{F}_2$ cuts} &\multicolumn{7}{c}{Adding $\mtc{F}_1$, $\mtc{F}_2$ and $\mtc{F}_3$ Cuts} \\
\hline
ID	&	$|S|$	&	$B$	&	CPU(s)	&	iters	&	obj	&	cutRelaxObj	&	solverCuts	&	$\mtc{F}_1$-cut	&	CPU(s)	&	iters	&	obj	&	cutRelaxObj	&	solverCuts	&	$\mtc{F}_1$-cut	&	$\mtc{F}_3$-cut	\\
\hline																															
FL1	&	10	&	5	&	4	&	5	&	4248	&	4313	&	21	&	29	&	5	&	5	&	4248	&	4312	&	21	&	29	&	11	\\
FL2	&	40	&	5	&	168	&	5	&	6499	&	6499	&	87	&	98	&	174	&	5	&	6499	&	6499	&	87	&	98	&	1	\\
FL3	&	40	&	10	&	11311	&	7	&	11923	&	12032	&	0	&	221	&	-	&	10	&	11926	&	12030	&	0	&	326	&	67	\\
FL4	&	40	&	15	&	3638	&	7	&	14249	&	14456	&	310	&	191	&	-	&	28	&	14250	&	14447	&	310	&	892	&	241	\\
FL5	&	60	&	5	&	234	&	5	&	7094	&	7094	&	96	&	92	&	271	&	5	&	7094	&	7094	&	96	&	92	&	0	\\
FL6	&	60	&	10	&	3695	&	6	&	13970	&	14034	&	231	&	196	&	-	&	20	&	13971	&	14033	&	231	&	757	&	78	\\
FL7	&	60	&	15	&	11582	&	7	&	19931	&	20166	&	556	&	322	&	-	&	9	&	19931	&	20163	&	556	&	441	&	103	\\
FL8	&	80	&	5	&	254	&	4	&	7409	&	7456	&	75	&	61	&	261	&	4	&	7409	&	7456	&	75	&	61	&	3	\\
FL9	&	80	&	10	&	4819	&	5	&	14404	&	14563	&	206	&	163	&	5088	&	5	&	14404	&	14561	&	206	&	165	&	23	\\
FL10	&	80	&	15	&	-	&	9	&	20281	&	20424	&	586	&	483	&	-	&	9	&	20281	&	20422	&	586	&	485	&	66	\\
FL11	&	100	&	5	&	1576	&	8	&	6905	&	6993	&	144	&	127	&	1411	&	8	&	6905	&	6991	&	144	&	127	&	17	\\
FL12	&	100	&	10	&	4452	&	6	&	14236	&	14369	&	211	&	197	&	-	&	17	&	14237	&	14349	&	211	&	640	&	112	\\
FL13	&	100	&	15	&	10125	&	6	&	20410	&	20485	&	519	&	323	&	10178	&	6	&	20410	&	20484	&	519	&	323	&	24	\\
FL14	&	200	&	5	&	2711	&	5	&	7377	&	7500	&	209	&	76	&	2586	&	5	&	7377	&	7498	&	209	&	76	&	17	\\
FL15	&	200	&	10	&	11244	&	7	&	14104	&	14339	&	298	&	220	&	-	&	9	&	14122	&	14322	&	298	&	306	&	50	\\
FL16	&	200	&	15	&	-	&	8	&	20850	&	21053	&	485	&	422	&	-	&	8	&	20850	&	21050	&	503	&	423	&	55	\\
FL17	&	300	&	5	&	12802	&	9	&	7532	&	7571	&	286	&	165	&	12847	&	9	&	7532	&	7571	&	286	&	165	&	19	\\
FL18	&	300	&	10	&	9011	&	5	&	15397	&	15483	&	277	&	177	&	9008	&	5	&	15397	&	15482	&	304	&	177	&	16	\\
FL19	&	300	&	15	&	12613	&	7	&	22776	&	22952	&	188	&	380	&	-	&	8	&	22776	&	22949	&	282	&	445	&	54	\\
FL20	&	400	&	5	&	-	&	8	&	7619	&	7697	&	237	&	136	&	-	&	8	&	7619	&	7695	&	239	&	136	&	14	\\
FL21	&	400	&	10	&	-	&	8	&	15445	&	15601	&	443	&	267	&	-	&	8	&	15448	&	15597	&	443	&	278	&	42	\\
FL22	&	400	&	15	&	-	&	8	&	22743	&	22985	&	552	&	421	&	-	&	8	&	22743	&	22982	&	552	&	431	&	66	\\
FL23	&	500	&	5	&	-	&	8	&	7674	&	7721	&	256	&	142	&	-	&	8	&	7674	&	7719	&	256	&	139	&	24	\\
FL24	&	500	&	10	&	-	&	8	&	15168	&	15489	&	385	&	231	&	-	&	8	&	15268	&	15483	&	385	&	243	&	66	\\
FL25	&	500	&	15	&	-	&	8	&	22766	&	22930	&	499	&	415	&	-	&	8	&	22703	&	22932	&	499	&	422	&	40	\\
FL26	&	1000	&	100	&	-	&	8	&	146474	&	147561	&	548	&	3056	&	-	&	8	&	148060	&	148821	&	640	&	3116	&	0	\\
FL27	&	1000	&	200	&	-	&	8	&	244054	&	246025	&	0	&	4057	&	-	&	8	&	244054	&	246025	&	0	&	4057	&	0	\\
FL28	&	1000	&	300	&	12605	&	7	&	317858	&	317911	&	830	&	5733	&	12604	&	7	&	317858	&	317910	&	830	&	5734	&	6	\\
FL29	&	1500	&	100	&	-	&	8	&	142000	&	143207	&	284	&	3118	&	-	&	8	&	143549	&	145698	&	284	&	3075	&	0	\\
FL30	&	1500	&	200	&	-	&	8	&	268452	&	272981	&	0	&	5144	&	-	&	8	&	264751	&	268994	&	0	&	5144	&	4	\\
\hline
\end{tabular}
}
\end{threeparttable}
\end{table}

\subsection{Practical Insights}
\label{sec:comp-manag-ins}
\begin{figure}[ht]
\centering
\subfigure[$\gamma_{ij}=0$, $\textrm{Obj}=3280.23$]{\includegraphics[width=0.35\textwidth]{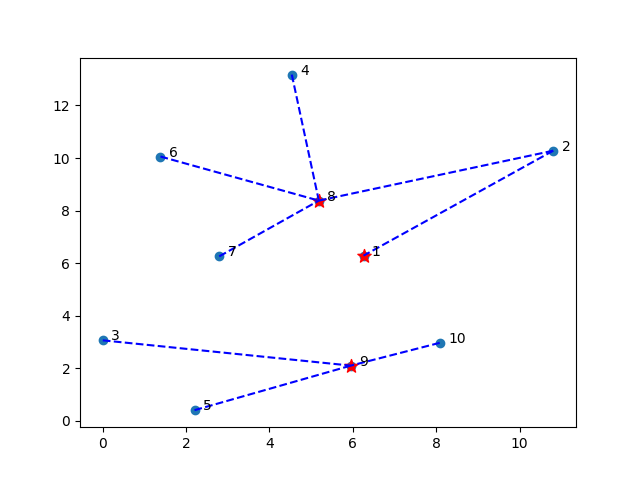}}
\subfigure[$\gamma_{ij}=0.2$, $\textrm{Obj}=3043.82$]{\includegraphics[width=0.35\textwidth]{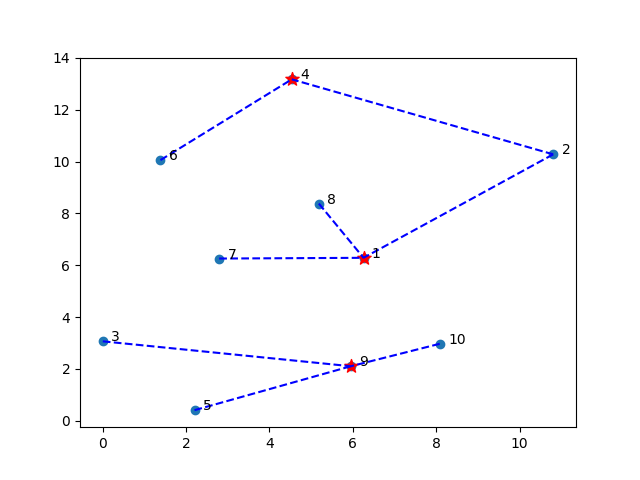}}
\subfigure[$\gamma_{ij}=0.4$, $\textrm{Obj}=2832.00$]{\includegraphics[width=0.35\textwidth]{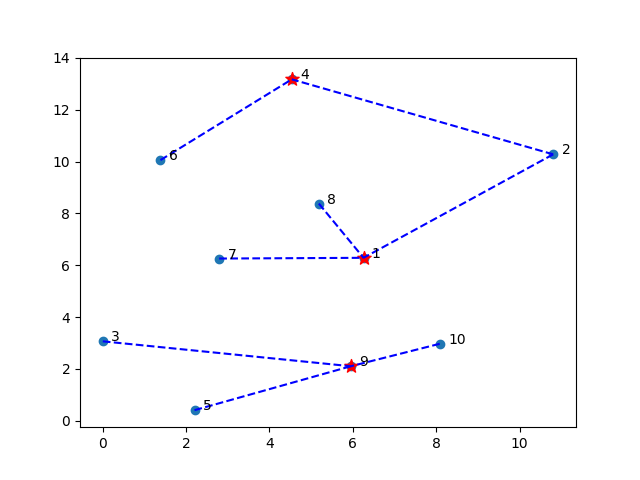}}
\subfigure[$\gamma_{ij}=0.8$, $\textrm{Obj}=2412.04$]{\includegraphics[width=0.35\textwidth]{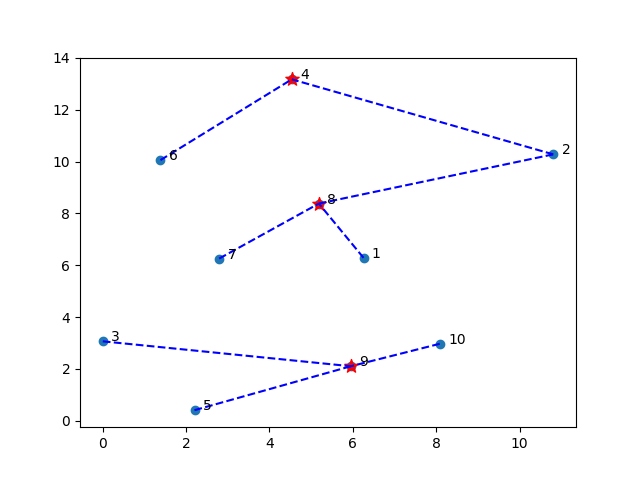}}
\caption{Impact of the level of ambiguity on the optimal service center location.  
All ten points (stars and circles) on the graph  represent the customer sites.
The points with stars are the optimal locations of service centers.  
The dashed lines indicate the customer flows from the customer sites to the service centers in the optimal solution.}
\label{fig:manag-insight}
\end{figure}

We now use a problem instance of \eqref{opt:RFL} to investigate the impact of the level of ambiguity 
on the optimal service center location decisions. The instance consists of 10 customer sites (denoted as L1-L10 respectively). 
All sites are candidate locations for the service centers. 
These sites are shown in Figure~\ref{fig:manag-insight}. 
The budget allows us to open 3 service centers. 
All parameters of this instance are created as described in Section~\ref{sec:gen-inst}. 
We tested the impact of the ambiguity level $\gamma^{ij}_2$ 
(see \eqref{eqn:sol-inner}) on the optimal service center location decisions.  
In the test we set $\gamma^{ij}_2=0, 0.2, 0.4, 0.6, 0.8$ for all $i\in S,j\in F$,
respectively to see how the robust optimal solution changes with the increment of the 
ambiguity level. Note that the setting $\gamma^{ij}_2=0$ is 
the nominal setting with no ambiguity in the utility function.

The robust optimal solutions of four settings $\gamma^{ij}_2=0, 0.2, 0.4, 0.8$ 
are shown in Figure~\ref{fig:manag-insight}. The plot for the setting $\gamma^{ij}_2=0.6$
is omitted since corresponding robust optimal solution is the same as the solution for setting $\gamma^{ij}_2=0.4$.
The optimal location decision in the nominal setting is at L1, L8 and L9.
For settings of $\gamma^{ij}_2=0.2, 0.4, 0.6$, 
the optimal location decision is L1, L4 and L9, which is different from the nominal setting.
When the parameter $\gamma^{ij}_2$ increases to 0.8, 
the optimal location decision changes to be L4, L8 and L9.

The example illustrates that by properly exploring parameter settings and solving \eqref{opt:RFL}, 
one may be able to identify the relevance of increasing precision to the data collection process.
In practice, if the robust optimal solution is very sensitive to the ambiguity level, 
it suggests the need for more accurately estimating the utility function before arriving at the final recommendation.

\section{A \cvd Testing Center Location Case Study}
\label{sec:covid19}
During a pandemic outbreak such as  \cvd screening tests are performed
to prevent a wider spread to decide on actions such as quarantine and contact tracing. Additionally, in such situations population needs to be rapidly vaccinated. We now discuss an application of \eqref{opt:RFL} to identifying locations where such tests/vaccinations can be performed. An optimal choice of locations 
should attract as many susceptible candidates as possible, and location convenience is a factor for at risk individuals.  As a case study, we consider location decisions for these test centers in San Diego county, CA. In the instance generation, we assume that each center has equal capacity and demand is generated synthetically using data available from the real-world as a basis.

\subsection{Parameters generation}
We use San Diego county zip codes as possible locations for setting up the testing/vaccination centers. The county has 78 zip codes. The distance between two zip code locations was calculated using
the latitude and longitude provided in \citep{zip-loc}. 
Locations with different zip codes are combined if their mutual distance is less than 0.05 miles. 
For illustrative purposes, we assume that a test/vaccination facility can be identified at the desired zip codes. 
After combining such zip codes  we have 69 locations that were used as model input.

The number of \cvd cases identified at each zip code was downloaded from \citep{case-san-diego}.  The 2010 census data provided resident populations by the zip code \citep{census}.
This data also provided 2010 average per capita income in each zip code. Using this data \cvd case rate per 10,000 residents was calculated based on the known number of cases in mid April 2020. For the purposes of this case study,  we scaled the total number of cases at each zip code proportionately to its population in our demand estimation. The model was developed to cover San Diego county population for testing/vaccination over a three month period. 
We set each testing center's capacity at 1,000 per day. For our model, 
we set the demand of each zip code to be 30\% of the population in that zip code. 

To generate the utility related parameters $\beta_{ij}$, $A_{ij}$ and $\Sigma_{ij}$
for all $i\in S, j\in F$ (the sets $S,F$ are the list of 69 zip code locations in this case), we simulated the survey conducted for each tuple $(i,j)$ and estimated these parameters 
using a linear regression model. We now describe how the input samples for the linear regression model were generated. Let $d_{ij}$ be the distance between locations $i$ and $j$, and $L_0$ be a distance threshold. The value of $L_0$ is set at the 20\% quantile of the sorted sequence $\Set*{d_{ij}}{i\in S,\;j\in F,\;i\neq j}$. Let us also define the relative income level at location $k$ as
\bdm
\textrm{IncRate}_k = \frac{\textrm{AvgInc}_k}{\max_{i\in F} \textrm{AvgInc}_i}, 
\edm
where the $\textrm{AvgInc}_k$ is the average per capita income at location $k$.
For every $i\in S$, we define the set $F_i=\Set*{j\in F}{d_{ij}\le L_0}$.
In the simulation experiments, for a tuple $(i,j)\in S\times F_i$
we generate $N$ samples of a random decision vector (a location profile) given as 
\beq\label{eqn:v-rand-loc}
\bs{v} = \bs{e}_j + \sum_{j^{\prime}\in F_i\setminus\{j\}} Z_{j^{\prime}}\bs{e}_{j^{\prime}},
\eeq
where $Z_{j^\prime}$ follows the Bernoulli distribution with $p=0.3$. 
It means that the decision vector locates a testing center at location $j$ and 
also locates a testing center at every $j^{\prime}\in F_i\setminus\{j\}$ 
randomly and independently with probability $0.3$.
We simulate the customer score for its location using 
\beq\label{eqn:u-rand-utility}
u_{ij}(\bs{v})=\sigma_1\max\Big\{0,1-\frac{d_{ij}}{L_0}\Big\}
+\sum_{j^{\prime}\in F_i\setminus\{j\}}Z_{j^{\prime}}\sigma_2\max\Big\{0,1-\frac{d_{ij^\prime}}{L_0}\Big\}
+\tau_1\textrm{IncRate}_j+\sum_{j^{\prime}\in F_i\setminus\{j\}} Z_{j^\prime}\tau_2\textrm{IncRate}_{j^\prime},
\eeq
where $Z_{j^\prime}$ is the same Bernoulli random parameter used in
specifying $\bs{v}$ in \eqref{eqn:v-rand-loc}. 
The parameters $\sigma_1,\sigma_2,\tau_1,\tau_2$
are truncated Gaussian random variables satisfying 
$\sigma_1\sim\max\{0,\;\mtc{N}(15,2)\}$, $\sigma_2\sim\max\{0,\;\mtc{N}(1,0.5)\}$, 
$\tau_1\sim\max\{0,\;\mtc{N}(3,1)\}$, $\tau_2\sim\max\{0,\;\mtc{N}(1,0.5)\}$.
Note that the first two terms in \eqref{eqn:u-rand-utility}
correspond to the contribution of distance to the utility
where shorter distance leads to higher utility contribution,
and the contribution of location $j^\prime$ will be zero if 
the distance $d_{ij^\prime}$ is larger than the 
threshold distance $L_0$ or $Z_{j^\prime}=0$.
The last two terms in \eqref{eqn:u-rand-utility}
correspond to the contribution of relative income level. The relationship in \eqref{eqn:u-rand-utility} is such that
the utility is positively correlated to the attractiveness 
of the location profile, where the attractiveness is measured
by the distance and income, i.e., individuals are more willing
to go for testing to a wealthy (safer) community that has a shorter distance to their home.
The location profile generated from \eqref{eqn:v-rand-loc}
and the utility \eqref{eqn:u-rand-utility} as the customer response to the location profile
are used to generate the data matrix and response vector in the linear 
regression model for estimating parameters $\beta_{ij}$, $A_{ij}$ and $\Sigma_{ij}$.
A pseudo-code for the procedure used in generating these parameters is given in Algorithm~\ref{alg:fit-linear-reg}.

\begin{algorithm}
	{\footnotesize
	\caption{\footnotesize An algorithm for simulating samples and fits a linear model to estimate 
	parameters $\beta_{ij}$, $A_{ij}$ and $\Sigma_{ij}$ for a specific tuple $(i,j)\in S\times F_i$.}
	\label{alg:fit-linear-reg}
	\begin{algorithmic}
	\State{Input: the number of samples $N$ and a tuple $(i,j)\in S\times F_i$.}
	\State{Set $V=\emptyset$ and $U=\emptyset$. Set $n\gets 0$.}
	\While{$n<N$}
		\State{Set $n\gets n+1$.}
		\State{Draw a sample $\bs{v}^n$ of location profile from \eqref{eqn:v-rand-loc} and 
		let $z^n_{j^\prime}$ be the realization of the Bernoulli variable $Z_{j^\prime}$ for $j^{\prime}\in F_i\setminus\{j\}$.}
		\State{Draw a sample of $\sigma_1,\sigma_2,\tau_1,\tau_2$ from the corresponding probability distributions, 
		respectively. Denote the realizations as $\sigma^n_1,\sigma^n_2,\tau^n_1,\tau^n_2$.}
		\State{Evaluate the utility using the formula \eqref{eqn:u-rand-utility} with the realization of the parameters
		$z^n_{j^\prime},\sigma^n_1,\sigma^n_2,\tau^n_1,\tau^n_2$. Denote the utility value as $u^n_{ij}$.}
		\State{Set $V\gets V\cup\{\bs{v}^n\}$, $U\gets U\cup\{u^n_{ij}\}$.}
	\EndWhile
	\State{Construct the data matrix $M=[\bs{v}^1,\ldots,\bs{v}^N]$ and 
		response vector $\bs{u}=[u^1_{ij},\ldots,u^N_{ij}]$.}
	\State{Fit the linear regression $M\beta_{ij}=\bs{u}+\bs{\varepsilon}$. 
	Let $\hat{\beta}_{ij}$ be the point estimation. Let $\Sigma_{ij}$ be
	the covariance matrix and let $A_{ij}$ be the diagonal matrix consisting of
	the diagonal elements of $\Sigma^{-1}_{ij}$. 
	Determine $b_{ij}$
	such that the ellipsoid $(\beta_{ij}-\hat{\beta}_{ij})^{\top}A_{ij}(\beta_{ij}-\hat{\beta}_{ij})\le b^2$ 
	is a 80\% confidence region of $\beta_{ij}$.}
	\State{Return $\hat{\beta}_{ij},A_{ij},\Sigma_{ij},b_{ij}$.}
	\end{algorithmic}
	}
\end{algorithm}

\subsection{Results analysis}
We created 36 instances of the \cvd testing center location problem. In all instances, the cost $b_j$ of locating a testing center at any location $j\in F$ is set to be one unit,
and the budget $B$ is the number of testing centers. 
We let $B$ range from 5 to 45 with an increment of 5, and let the options for the number of samples $N$
(the input of Algorithm~\ref{alg:fit-linear-reg}) be 500, 1000, 1500 and 2000. 
The 36 numerical instances are created corresponding to 36 combinations of
$(B,N)$ parameters. Note that Algorithm~\ref{alg:fit-linear-reg} has been 
run for every tuple $(i,j)\in S\times F_i$ to generate $N$ samples for linear regression
in generating an instance. Each numerical instance is solved using two methods: 
the lifted MISOCP formulation and the cutting-plane method with a 4-hour CPU time limit.
The numerical results are compared and summarized in Table~\ref{tab:covid-19}.
All instances are solved to optimality by the lifted MISOCP formulation within the time limit.
It is observed that the solution time of using the cutting-plane method is magnitude smaller
than the lifted formulation while it returns the same objective value as the lifted formulation
in 31 instances.

\begin{figure}
\centering
\includegraphics[scale=0.35]{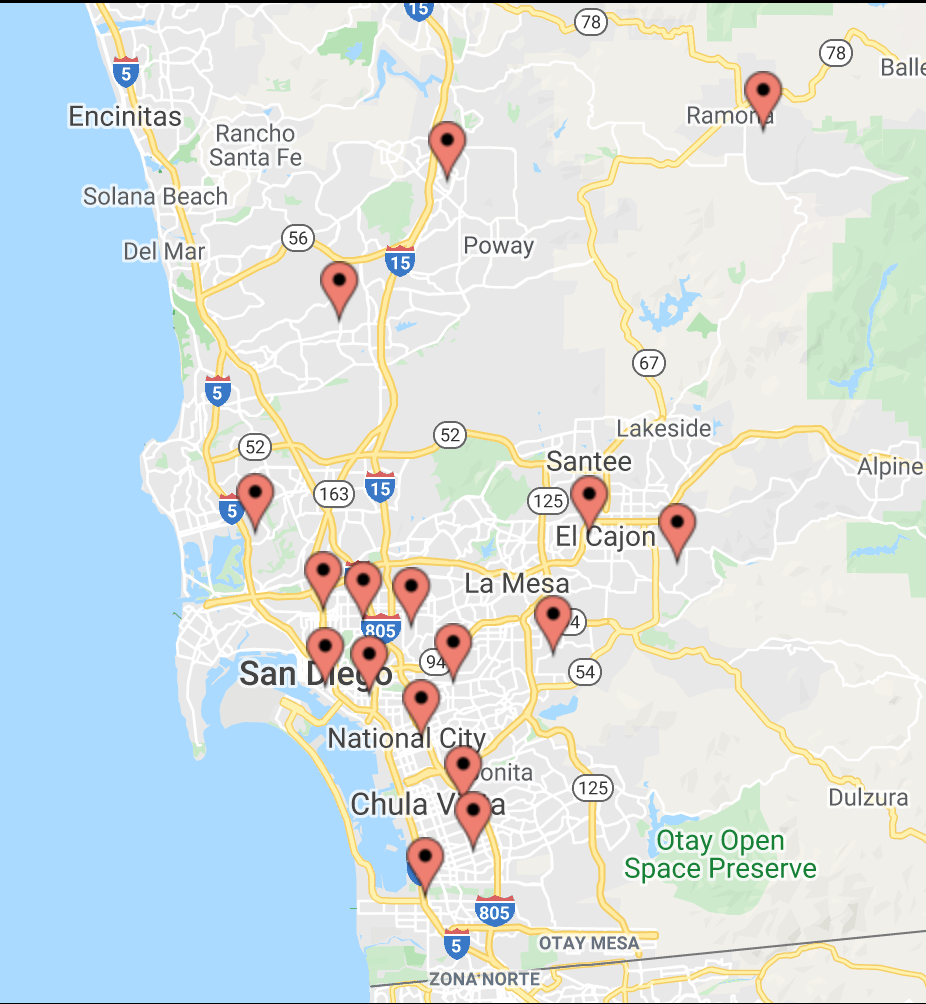}
\caption{The optimal locations of 20 \cvd testing centers for the case that $N=1000$ and $B=20$.}
\label{fig:location-plot-on-map}
\end{figure}

As an illustration, the optimal locations (zip codes) of the \cvd testing centers for the instance $N=1000, B=20$
are shown in Figure~\ref{fig:location-plot-on-map}. In this optimal solution, 10  centers are located in the City of San Diego, 4 centers are located in the nearby suburbs, and 3 centers are located
in different towns to the north of the city. 

\begin{figure}
\centering
\includegraphics[width=0.45\textwidth]{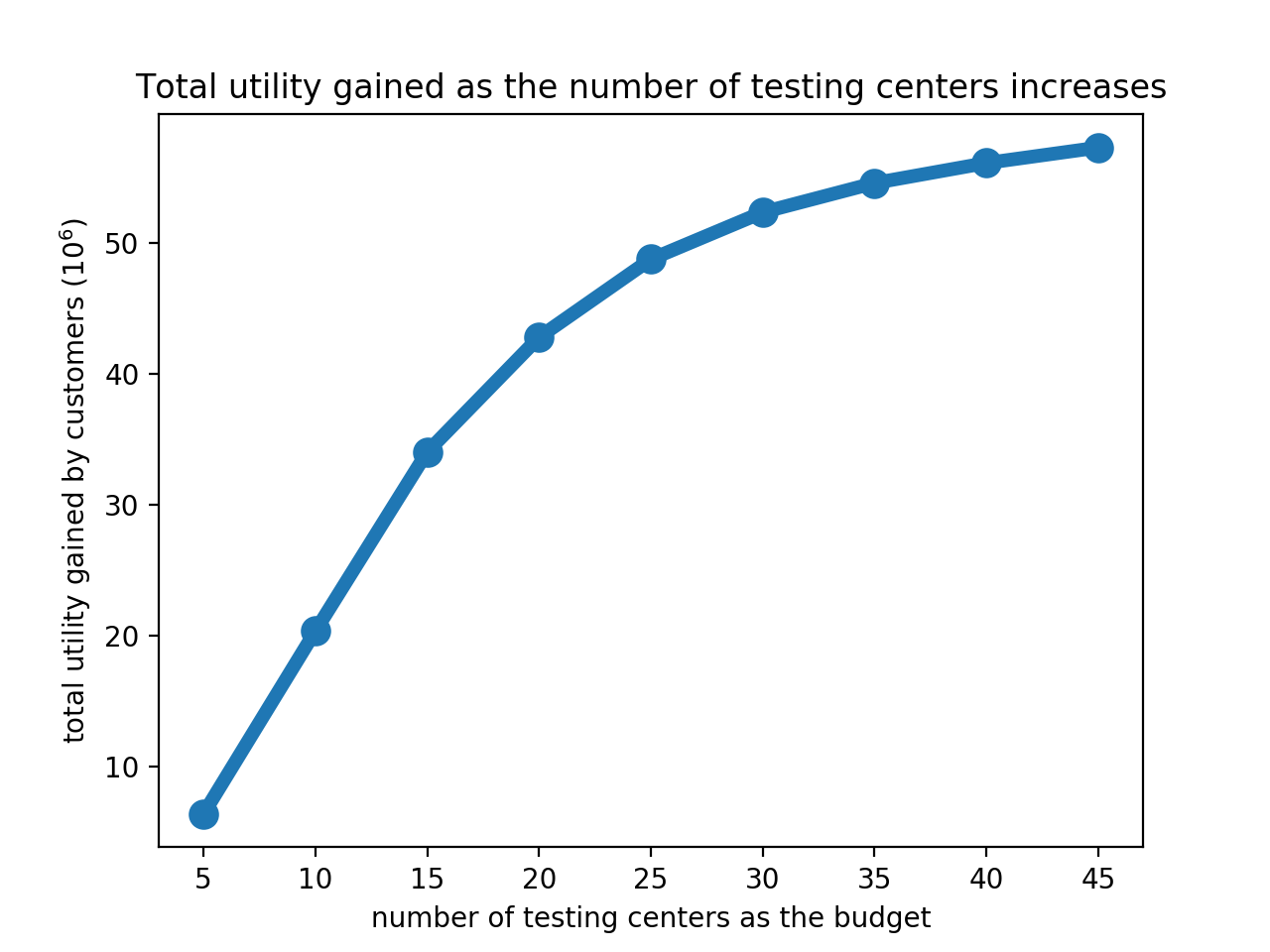}
\caption{Total utility with $N=1000$ versus $B$.}\label{fig:obj-vs-B}
\end{figure}

We now discuss the relationship between the number of  facilities $B$ and their maximum utility.
Figure~\ref{fig:obj-vs-B} shows the optimal  utility with increasing value of $B$. It is observed that in the range $1\le B\le 15$, 
the objective value increase almost linearly with the number of centers. In the $B\ge 15$ range, the objective value still increases with $B$, but with a slower rate of increase. This dependency behaves as a 
concave function defined on discrete points (the number of locations). For
$B\le 15$, the capacity of all centers is a bottleneck, resulting in the linear behavior. For $B\ge 20$, the  center capacity is no longer a bottleneck. It implies that the utility gains from adding a new center diminish as some of the centers are not fully utilized regardless of their location. 

We also studied the impact of sample size $N$ used for fitting the utility model on the objective value.
For a given budget, we report the objective values corresponding to four different sample size in Table~\ref{tab:obj-B-N}. We observe that between
$N=500$ and $N=1000$ the change in the objective value  is small (0.41\% change on average), in comparison to the change for  $N=50$ and $N=100$ (2.4\% change on average). It indicates the convergence in the objective value of the problem with an increase in sample size as tighter confidence intervals are now available resulting in a reduced ambiguity set.

\begin{table}
		{\centering
		\caption{ Total utility versus the budget for 4 options of sample size. 
		Obj1, Obj2, Obj3 and Obj4 are total utility corresponding to sample size 50, 100, 500 and 1000, respectively.}
		\label{tab:obj-B-N}
		\centering
		\begin{tabular}{c|rrrrrrrrr}
		\hline\hline
		Budget($B$)	&	5	&	10	&	15	&	20	&	25	&	30	&	35	&	40	&	45	\\
		\hline
Obj1	&	6.28	&	20.11	&	30.11	&	41.58	&	44.28	&	48.34	&	51.94	&	51.94	&	55.27	\\
Obj2	&	6.21	&	20.15	&	31.46	&	41.63	&	45.39	&	52.13	&	52.24	&	53.86	&	56.08	\\
Obj3	&	6.35	&	20.08	&	33.81	&	42.98	&	48.85	&	52.25	&	54.45	&	55.99	&	57.11	\\
Obj4	&	6.34	&	20.10	&	34.02	&	42.79	&	48.76	&	52.35	&	54.59	&	56.12	&	57.30	\\
		\hline
		\end{tabular}
		}
\end{table}

\begin{table}
\centering
{\scriptsize
 \caption{\scriptsize Computational results of numerical instances for \cvd testing center location in the San Diego county.
 Each instance is under solving for 4 hours with two approaches: the lifted MISOCP formulation and the cutting-plane method. All instances have been solved to optimality using the lifted MISOCP formulation within 4 hours.  
Since the cutting-plane method cannot find an optimal solution in general, the best objective value is reported. 
For the numerical instances with $(|S|, B, N)=(69, 15, 500), (69, 20, 500), (69, 40, 500)$, Gurobi encountered 
numerical troubles when solving the projection problem \eqref{opt:min-dist}, and hence it is not able to return a qualified
solution of the instance. For these cases, we put `-' at each column for the cutting-plane method as a place holder.
}
 \label{tab:covid-19}
 \begin{tabular}{ccc|cccc|cccccc}
 \hline\hline
 & & & \multicolumn{4}{c|}{Lifted Formulation} & \multicolumn{6}{c}{Cutting-Plane Method} \\
 \hline
$|S|$	&	$B$	&	Samples	&	solTime(s)	&	obj	&	nodes	&	solverCuts	&	solTime(s)	&	iters	&	obj	&	nodes	&	solverCuts	&	devCuts	\\
\hline
69	&	5	&	500	&	9585	&	273881	&	7095	&	6	&	65	&	2	&	273881	&	121	&	4	&	5	\\
69	&	5	&	1000	&	9638	&	282014	&	8687	&	13	&	120	&	2	&	282014	&	118	&	1	&	5	\\
69	&	5	&	1500	&	7991	&	286858	&	6723	&	17	&	108	&	2	&	286858	&	104	&	3	&	5	\\
69	&	5	&	2000	&	8392	&	281898	&	8638	&	12	&	222	&	4	&	281898	&	119	&	5	&	15	\\
\hline
69	&	10	&	500	&	12205	&	755400	&	8295	&	6	&	411	&	6	&	751288	&	148	&	8	&	50	\\
69	&	10	&	1000	&	10637	&	765988	&	5733	&	0	&	1639	&	6	&	765988	&	108	&	3	&	53	\\
69	&	10	&	1500	&	11382	&	764279	&	7363	&	0	&	1857	&	7	&	764279	&	123	&	3	&	63	\\
69	&	10	&	2000	&	8831	&	767109	&	6444	&	11	&	1383	&	6	&	767109	&	114	&	4	&	53	\\
\hline
69	&	15	&	500	&	7631	&	1285211	&	2287	&	6	&	-	&	-	&	-	&	-	&	-	&	-	\\
69	&	15	&	1000	&	1449	&	1304416	&	1305	&	23	&	1801	&	4	&	1304416	&	273	&	0	&	48	\\
69	&	15	&	1500	&	1366	&	1316221	&	1220	&	24	&	1426	&	4	&	1316221	&	161	&	0	&	48	\\
69	&	15	&	2000	&	5473	&	1314293	&	1915	&	12	&	430	&	2	&	1314293	&	198	&	0	&	16	\\
\hline
69	&	20	&	500	&	8523	&	1738786	&	5551	&	6	&	-	&	-	&	-	&	-	&	-	&	-	\\
69	&	20	&	1000	&	2400	&	1759829	&	2272	&	20	&	309	&	2	&	1759829	&	897	&	11	&	21	\\
69	&	20	&	1500	&	2188	&	1771388	&	1999	&	26	&	955	&	4	&	1771388	&	618	&	10	&	63	\\
69	&	20	&	2000	&	2099	&	1767474	&	2232	&	33	&	976	&	4	&	1767474	&	822	&	25	&	63	\\
\hline
69	&	25	&	500	&	10417	&	2287932	&	5189	&	5	&	860	&	5	&	2287932	&	910	&	10	&	100	\\
69	&	25	&	1000	&	9698	&	2300836	&	4450	&	52	&	586	&	4	&	2300836	&	1127	&	21	&	75	\\
69	&	25	&	1500	&	9811	&	2304030	&	4164	&	53	&	799	&	5	&	2304030	&	666	&	13	&	100	\\
69	&	25	&	2000	&	7454	&	2300852	&	3318	&	49	&	544	&	4	&	2300852	&	1550	&	10	&	75	\\
\hline
69	&	30	&	500	&	221	&	2909555	&	793	&	39	&	66	&	4	&	2909556	&	100	&	3	&	90	\\
69	&	30	&	1000	&	272	&	2930043	&	1035	&	80	&	56	&	4	&	2930043	&	71	&	6	&	90	\\
69	&	30	&	1500	&	320	&	2944779	&	1402	&	92	&	34	&	2	&	2944779	&	134	&	6	&	30	\\
69	&	30	&	2000	&	298	&	2937953	&	1238	&	39	&	38	&	2	&	2937953	&	87	&	14	&	30	\\
\hline
69	&	35	&	500	&	844	&	3433722	&	1891	&	85	&	154	&	6	&	3433722	&	144	&	26	&	175	\\
69	&	35	&	1000	&	5062	&	3462695	&	2396	&	15	&	95	&	4	&	3462695	&	131	&	3	&	105	\\
69	&	35	&	1500	&	4536	&	3471085	&	1626	&	50	&	114	&	5	&	3471085	&	85	&	4	&	140	\\
69	&	35	&	2000	&	3545	&	3472486	&	1748	&	25	&	102	&	4	&	3472489	&	178	&	6	&	105	\\
\hline
69	&	40	&	500	&	182	&	4028719	&	26	&	38	&	-	&	-	&	-	&	-	&	-	&	-	\\
69	&	40	&	1000	&	169	&	4045728	&	55	&	99	&	21	&	2	&	4045728	&	1	&	1	&	40	\\
69	&	40	&	1500	&	279	&	4054252	&	61	&	79	&	17	&	2	&	4054252	&	1	&	6	&	40	\\
69	&	40	&	2000	&	254	&	4058904	&	31	&	57	&	20	&	2	&	4058904	&	1	&	3	&	40	\\
\hline
69	&	45	&	500	&	104	&	4497358	&	44	&	25	&	11	&	2	&	4497358	&	1	&	3	&	45	\\
69	&	45	&	1000	&	170	&	4512684	&	92	&	27	&	68	&	5	&	4512684	&	1	&	3	&	180	\\
69	&	45	&	1500	&	177	&	4522907	&	157	&	21	&	79	&	6	&	4522907	&	1	&	7	&	225	\\
69	&	45	&	2000	&	76	&	4523842	&	52	&	24	&	55	&	5	&	4523904	&	1	&	2	&	180	\\
\hline
 \end{tabular}
}
\end{table}

\section{Concluding Remarks}
The utility-robust facility location model captures the 
endogenous uncertainty of customers' utility in decision making.
The moment-based ambiguity set constructed in this paper 
for the decision dependent utility leads to a mixed 0-1 second-order-cone program.
This reformulation shows that the discrete optimization models with decision dependent ambiguity sets may admit 
a convex reformulation with mixed-binary variables. 
Incorporating the convexification cuts developed in this paper helps solve the \eqref{opt:RFL} problem more efficiently, especially
for large instances where the approach without the cuts can not achieve desired four digit accuracy in the solution within the time limit of four hours. In practice, the ambiguity level can be determined empirically 
based on estimation accuracy of the parameters from the collected data. Moreover, by using several values of this parameter we can test the sensitivity of the optimal solution. The illustrative example of locating Covid-19 centers in San Diego county, CA reveals that the optimal objective value of the model is concave at discrete values of the number of test centers. This example also confirms the effectiveness of adding identified cuts in closing the optimality gap and generating an improved solution when the computational time budget is limited.

This paper assumed that the utility function is linear, 
and a moment based model for describing the ambiguity set for decision dependent utilities. Alternative models for expressing a decision maker's utility may be explored in the future. 

Stochastic optimization framework to model uncertain demand has been proposed for the facility location problems \citep{2006-fac-loc-uncert-rev}.   We now present a generalization of the basic \eqref{opt:RFL} model for the case where
the customer demand is stochastic with a finite support. In this case, the customer demand is denoted as $D^{\omega}_i$ to represent the demand value at
scenario $\omega\in\Omega$, and the number of customers going to a facility is denoted by $x^{\omega}_{ij}$. In the stochastic demand case we can further define an ambiguity set $\mP^{\Omega}$ for the unknown probability distribution over scenarios. With one more layer of ambiguity on the probability distributions over scenarios,  
the model \eqref{opt:RFL} is then formulated as a distributionally-robust two-stage stochastic optimization
problem written as follows:
\begin{equation}\label{opt:SD-RFL}
\begin{aligned}
&\underset{\bs{y}}{\textrm{max}}\;\;\bs{h}^{\top}\bs{y}+\underset{P\in\mathcal{P}^{\Omega}}{\textrm{min}}\E_{P}[\mathcal{Q}(\bs{y},\omega)] \\
&\textrm{ s.t. } \sum_{j\in F} b_j y_j \le B,   \quad y_j\in\{0,1\}\;\forall j\in F,
\end{aligned}
\tag{SD-RFL}
\end{equation}
where the recourse function $\mathcal{Q}(\bs{y},\omega)$ is for the scenario $\omega\in\Omega$,
and it is defined similarly as \eqref{opt:RFL-II} with a scenario index on the demand and the number of customers going to a facility for service.
As discussed in Section~\ref{sec:mixed-bin-conic}, the \eqref{opt:RFL} model admits 
a mixed 0-1 second-order-cone program (MISOCP) reformulation based on a definition of $\Puyij$ using moments.
Similarly, the \eqref{opt:SD-RFL} model can be reformulated as a distributionally-robust two-stage stochastic
mixed 0-1 second-order-cone program (DR-TSS-MISOCP). Solving such problems effectively requires further algorithmic development. A decomposition branch-and-bound method for solving a general 
DR-TSS-MISOCP problem is developed in our recent work \citep{luo2019-BB-two-stage-0-1-socp}. This algorithm is used to solve small instances of \eqref{opt:SD-RFL} in \citep{luo2019-BB-two-stage-0-1-socp}.
The reformulations and convexification cuts
developed in the current paper were used to strengthen the second-stage problem of \eqref{opt:SD-RFL} in the numerical 
study conducted in \citep{luo2019-BB-two-stage-0-1-socp} with a significant improvement in the computational performance. We refer the reader to our companion paper \citep{luo2019-BB-two-stage-0-1-socp} for a more detailed discussion on this topic. 

\if0\blind{
\section*{Acknowledgements}
This research was partially supported by the ONR grant N00014-18-1-2097. The author would like
to thank the anonymous referee and the associate editor who provided useful and constructive comments that lead
to a solid improvement of the manuscript.} \fi

\bibliographystyle{chicago}
\spacingset{1}
\bibliography{reference-database-current}
	
\end{document}